\documentclass[12pt,a4paper]{amsart}

\usepackage{amssymb}
\usepackage{amsmath}
\usepackage{mathrsfs} % needed for the curly 'F' and 'E' and 'O'(with \mathscr)
\usepackage{pst-all}  % to make the graphic
\usepackage[dvips,final]{graphicx}

\theoremstyle{plain}
\newtheorem{theorem}{Theorem}

\newtheorem{corollary}[theorem]{Corollary}
\newtheorem{lemma}[theorem]{Lemma}

\theoremstyle{definition}
\newtheorem{definition}[theorem]{Definition}
\newtheorem{example}[theorem]{Example}
\newtheorem{remark}[theorem]{Remark}

\newcommand{\bm}{\overline{M}}
\newcommand{\bg}{\overline{g}}
\newcommand{\bnabla}{\overline{\nabla}}

\newcommand{\dd}{\mathrm{d}}

\newcommand{\E}{\mathscr{E}}

\newcommand{\T}{\mathrm{T}}
\newcommand{\Order}{\mathscr{O}}
\newcommand{\dA}{\mathrm{det}A}
\newcommand{\area}{\mathrm{Area}}
\newcommand{\length}{\mathrm{Length}}

\newcommand{\II}{\mathrm{I\!I}}
\newcommand{\III}{\mathrm{I\!I\!I}}

\newcommand{\Hess}{\mathrm{Hess}}
\newcommand{\R}{\overline{\mathrm{R}}}
\newcommand{\Ric}{\overline{\mathrm{Ric}}}

\newcommand{\deltasnought}{\left.\frac{\partial}{\partial s}\right\vert_{s=0}}

\newcommand{\tensorheight}{1.2mm}

\newcommand{\ina}{\overline{\rule{0pt}{\tensorheight}a}}
\newcommand{\inc}{\overline{\rule{0pt}{\tensorheight}c}}
\newcommand{\ine}{\overline{\rule{0pt}{\tensorheight}e}}
\newcommand{\ini}{\overline{\rule{0pt}{\tensorheight}\imath}}
\newcommand{\inj}{\overline{\rule{0pt}{\tensorheight}\jmath}}
\newcommand{\inm}{\overline{\rule{0pt}{\tensorheight}m}}
\newcommand{\ins}{\overline{\rule{0pt}{\tensorheight}s}}
\newcommand{\inu}{\overline{\rule{0pt}{\tensorheight}u}}
\newcommand{\inv}{\overline{\rule{0pt}{\tensorheight}v}}
\newcommand{\inw}{\overline{\rule{0pt}{\tensorheight}w}}

\font \smallindexfont = cmbx5
\newcommand{\innul}{\overline{\rule{0pt}{\tensorheight}\textrm{\smallindexfont 0}}}
\newcommand{\ineen}{\overline{\rule{0pt}{\tensorheight}\textrm{\smallindexfont 1}}}

\begin{document}

\title[The Mean Curvature of the Second Fundamental Form]{The Mean Curvature\\ of the Second Fundamental Form\\ of a Hypersurface}

\author[Stefan Haesen, Steven Verpoort]{Stefan Haesen$^{\star}$, Steven Verpoort$^{\bullet}$}

\thanks{$^{\star}$ S. Haesen was a postdoctoral researcher at the K.U.Leuven while the work was initiated. He was partially supported by the Spanish MEC Grant MTM2007-60731 with FEDER funds and the Junta de Andaluc\'\i a Regional Grant P06-FQM-01951.}

\thanks{$^{\bullet}$ Corresponding author.}

\thanks{$^{\star}{}^{\bullet}$ 
Both authors were partially supported by the Research Foundation -- Flanders (project G.0432.07).}


%%%%%%%%%%%% Authors' addresses %%%%%%%%%%%%%
\address{% First author
\textit{Stefan Haesen:} \endgraf
Simon Stevin Institute for Geometry \endgraf 
Wilhelminaweg 1 \endgraf
nn 2042 Zandvoort \endgraf
The Netherlands. \endgraf
}
\email{stefan.haesen@geometryinstitute.org}

\address{% Second author
\textit{Steven Verpoort:} \endgraf
K.U.Leuven \endgraf 
Departement Wiskunde \endgraf
Afdeling Meetkunde \endgraf
Celestijnenlaan 200B bus 2400 \endgraf
3001 Heverlee \endgraf
Belgium. \endgraf
}
\email{steven.verpoort@gmail.com}

%%%%%%%%%%%%%%%%%%%%%%%%%%%%%%%%%%%%%%%%%%%%%%%%%%%%%%%%%%%%%%%%%%%%%%%%%%%%%%%%%%%%%%%%%%%%%
\begin{abstract}
An expression for the first variation of the area
functional of the second fundamental form is given for a hypersurface in a
se\-mi-Rie\-mann\-ian space. The concept of
the ``\textit{mean curvature of the second fundamental form}'' is then introduced. Some
characterisations of extrinsic hyperspheres in terms of this curvature are
given.

\vspace{3mm}
\noindent\textsc{Notice.} This article (ArXiv:0709.2107v3) is an extended version of \cite{thisarticle_shortversion}.
\end{abstract}
%%%%%%%%%%%%%%%%%%%%%%%%%%%%%%%%%%%%%%%%%%%%%%%%%%%%%%%%%%%%%%%%%%%%%%%%%%%%%%%%%%%%%%%%%%%%%

\maketitle

%%%%%%%%%%%%%%%%%%%%%%%%%%%%%%%%%%%%%%%%%%%%%%%%%%%%%%%%%%%%%%%%%%%%%%%%%%%%%%%%%%%%%%%%%%%%%

\section{Introduction and Outline of the Article}

We shall be concerned with hypersurfaces of a se\-mi-Rie\-mann\-ian manifold,
for which the real-valued second fundamental form $\II$ is a
se\-mi-Rie\-mann\-ian metrical tensor. 
The geometry of such hypersurfaces can be explored with respect to either the first or the second
fundamental form. 

In analogy with the classical study of the geometry of hypersurfaces as determined by their first fundamental form, a  distinction can be made between the \textit{intrinsic geometry of the second
fundamental form}, which is determined by measurements of $\II$-lengths on the hypersurface only,
and the \textit{extrinsic geometry of the second fundamental form}, which concerns those measurements for which the geometry of the
second fundamental form of the hypersurface is compared with the corresponding geometry of nearby hypersurfaces.

\subsection*{The Intrinsic Geometry of the Second Fundamental Form.}$\ $

\noindent
It is a natural question to investigate the relation between the \textit{intrinsic geometry
of the second fundamental form} and the shape of the original hypersurface,
and for this purpose the intrinsic curvatures of the second fundamental form have already been studied.

For example, it is well-known that the second fundamental form is a flat Lorentzian metric on a minimal surface in $\mathbb{E}^3$.
Conversely, D.E. Blair and T. Koufogiorgos showed that the pieces of a helicoid are the only non-developable ruled surfaces in $\mathbb{E}^3$ for which the Gaussian curvature $K_{\II}$ of $(M,\II)$ vanishes \cite{blair_kouf}, and M. Becker and W. K\"uhnel showed that the catenoid is the only surface of revolution in $\mathbb{E}^3$ for which the second fundamental form is a complete and flat Lorentzian metric \cite{beckerkuhnel}.

Numerous related characterisations of Euclidean spheres among ovaloids, i.e., among compact hypersurfaces in a Euclidean space with a positive definite second fundamental form, have already been found.
For example, R. Schneider's theorem characterises the hyperspheres as the only ovaloids for which the second fundamental form has constant sectional curvature \cite{schneider1972}. Some generalisations of this theorem for surfaces in certain Lorentzian manifolds have been found by J.A. Aledo, A. Romero, \textit{et al}
\cite{aledoromero2003,aledohaesenromero}.

However, in the present article we are 
not concerned with the geometry of the second fundamental form
from the intrinsic point of view, but we will study an aspect of the 
\textit{extrinsic geometry of the second fundamental form}.

\subsection*{The Extrinsic Geometry of the Second Fundamental Form.}$\ $

\noindent
As is known, the mean curvature $H$ of a hypersurface of a se\-mi-Rie\-mann\-ian manifold describes the instantaneous response of the area functional 
with respect to deformations of the hypersurface. Since we are studying hypersurfaces for which the second fundamental form is a se\-mi-Rie\-mann\-ian metrical
tensor, areas can be measured with respect to the second fundamental form as well, so we can associate
to any such hypersurface $M$ its area as
measured in the geometry of the second fundamental form.
This area, which will be denoted by $\area_{\II}(M)$, is related
to the classical area element $\dd\Omega$ by
\[
\area_{\II}(M)=\int_M \sqrt{|\dA|}\,\dd\Omega\,,
\]
where $A$ denotes the shape operator of the hypersurface.

In this article, the notion of mean curvature will be tailored
to the geometry of the second fundamental form: the function which measures the rate of change of $\area_{\II}(M)$ under a deformation of $M$, will be called the
\textit{mean curvature of the second fundamental form} and denoted by $H_{\II}$. In this way, a concept which belongs to the \textit{extrinsic geometry of the second fundamental form} will be introduced
in analogy with a well-known concept in the classical theory of hypersurfaces.
The mean curvature of the second fundamental form was defined originally by E. Gl\"assner \cite{glassnersimon,glassner} for surfaces in $\mathbb{E}^3$.
The corresponding variational problem has been studied by F. Dillen and W. Sodsiri \cite{dillensodsiri2005} for surfaces in $\mathbb{E}^{3}_1$, and for Riemannian surfaces in a three-dimensional se\-mi-Rie\-mann\-ian manifold in
\cite{haesenverpoortverstraelen}.

Some characterisations of the spheres in which this curvature $H_{\II}$ is involved have been found. For example,
it has been shown that the spheres are the only ovaloids in $\mathbb{E}^3$ which satisfy $H_{\II}=C \sqrt{K}$; furthermore,
the spheres are the only ovaloids on which $H_{\II}-K_{\II}$ does not change sign (see \cite{verpoort} and G. Stamou's \cite{stamou2003}).

In the initiating \S\,\ref{sec:def} of this article, the notation will be explained and several useful formulae from the theory of hypersurfaces will be briefly recalled.

In the following \S\,\ref{sec:IImin}, the first variation of the area functional of the second fundamental form is calculated and the mean curvature of the second fundamental form is defined. 

In \S\,\ref{sec:somefirstresults}, a comparison result for the Levi-Civita connections of the first and the
second fundamental form, which will be used in some of the subsequent proofs, is established.

In the subsequent sections (\S\S\,\ref{sec:spaceform}--\ref{sec:3d})
the mean curvature of the second fundamental form will be employed to give 
several characterisations of extrinsic hyperspheres as 
the only hypersurfaces in space forms, an Einstein space, and a three-dimensional manifold, respectively,
which can satisfy certain inequalities in which the mean curvature
of the second fundamental form is involved.

In \S\,\ref{sec:curves} the expresssion for $H_{\II}$
will be investigated for curves. This is of particular interest, since the length of the second fundamental form
of a curve $\gamma$,
\[
\length_{\II}(\gamma)=\int \sqrt{|\kappa|}\,\dd s\,,
\]
(where $\kappa$ is the geodesic curvature and $s$ an arc-length parameter) is a modification of the
classical bending energy
\[
\int \kappa^{2}\,\dd s
\]
which has already been studied by D. Bernoulli and L. Euler. Moreover, the results we present agree with W. Blaschke's description of J. Radon's variational problem \cite{blaschke} and with a more recent article of
J. Arroyo, O.J. Garay and J.J. Menc\'\i a \cite{arroyogaraymencia2003}.

In the final section, \S\,\ref{sec:geodhyp},
the function $H_{\II}$ will be investigated
for (sufficiently small) geodesic hyperspheres in a Riemannian manifold
by means of the method of power series expansions, which was applied extensively by A. Gray \cite{gray1973},
and also by B.-Y. Chen and L. Vanhecke \cite{chenvanhecke,grayvanhecke}.
Furthermore, we address the question of whether the locally flat spaces are characterised by the property that every geodesic hypersphere has the same $\II$-area as a Euclidean hypersphere with the same radius. 

%%%%%%%%%%%%%%%%%%%%%%%%%%%%%%%%%%%%%%%%%%%%%%%%%%%%%%%%%%%%%%%%%%%%%%%%%%%%%%%%%%%%%%%%%%%%%

\section{Definitions, Notation, and Useful Formulae}
\label{sec:def}

\subsection{Assumption.} All hypersurfaces are understood to be embedded and connected.

\subsection{Nomenclature.}
A hypersurface in a se\-mi-Rie\-mann\-ian manifold is said to be \makebox{(semi-)}\hspace{0mm}Riemannian if the
restriction of the metric to the hypersurface is a \makebox{(semi-)}Rie\-mann\-ian metrical tensor.

\subsection{Notation.}
\label{remarkbarIInotation}
Since a hypersurface $M$ in a manifold
$\bm$ will be studied, geometric objects in $\bm$ are
distinguished from their analogues in $M$ with a
$\overline{\mathrm{bar}}$. Geometric entities derived from the
second fundamental form are distinguished from those derived from
the first fundamental form by means of a sub- or superscript ${\II}$.
For example, the area element obtained from the second fundamental form will be written as
$\dd \Omega_{\II}$.

\subsection{Notation.}
The set of all vector fields on a manifold $M$ will be
denoted by $\mathfrak{X}(M)$. Furthermore, $\mathfrak{F}(M)$ stands
for the set of all real-valued functions on $M$. If $(M,g)$ is a se\-mi-Rie\-mann\-ian
submanifold of a se\-mi-Rie\-mann\-ian manifold $(\bm,\bg)$, the set of all vector fields on $M$ which
take values in the tangent bundle $\mathrm{T}\bm$ is denoted by
$\overline{\mathfrak{X}}(M)$. The orthogonal projection $\mathrm{T}_p\overline{M}\rightarrow
\mathrm{T}_p M$ will be denoted by $\left[\cdot\right]^{T}$.

\subsection{The Laplacian.}
The sign of the Laplacian will  be chosen so that $\Delta f=f''$ for a real-valued function on $\mathbb{R}$.

\subsection{The fundamental forms.}
Let $M$ be a se\-mi-Rie\-mann\-ian hypersurface of dimension $m$ in a se\-mi-Rie\-mann\-ian manifold
$(\bm,\bg)$. We shall
suppose that a unit normal vector field
$U\in\overline{\mathfrak{X}}(M)$ has been chosen  on $M$. The shape operator $A$, the second
fundamental form $\II$ and the third fundamental form $\III$
of the hypersurface $M$ are defined by the formulae
\begin{equation}
\left\{
\begin{array}{ccccccccc}
A    & : & \mathfrak{X}(M)                       & \rightarrow & \mathfrak{X}(M) & : & V    & \mapsto & -\bnabla_{V}U   \,;\\
\II  & : & \mathfrak{X}(M) \times\mathfrak{X}(M) & \rightarrow & \mathfrak{F}(M) & : &(V,W) & \mapsto & \alpha\,g(A(V),W) \,;\\
\III & : & \mathfrak{X}(M)\times\mathfrak{X}(M)  & \rightarrow & \mathfrak{F}(M) & : &(V,W) & \mapsto & g(A(V),A(W))      \,,
\end{array}
\right. \label{eq:defII}
\end{equation}
where $\alpha=\bg( U,U)=\pm1$.
It will be assumed that the second fundamental form is a se\-mi-Rie\-mann\-ian metric on $M$.

\subsection{Frame fields.}\label{sec:frame_fields}
Let $\{E_1,\ldots,E_m\}$ denote a frame field on $M$ which is orthonormal with respect to
the first fundamental form $g$. Define $\varepsilon_i$ ($i=1,\ldots,m$) by $\varepsilon_i=g(E_i,E_i)=\pm1$.
Furthermore, let $\{V_1,\ldots,V_m\}$ be a frame field on $M$ which is  orthonormal with respect to
the second fundamental form $\II$. Define $\kappa_i$ ($i=1,\ldots,m$) by $\kappa_i=\II(V_i,V_i)=\pm1$.

\subsection{Curvature.}
The following convention concerning the Riemann-Christoffel
curvature tensor $\textrm{R}$ will be made: for
$X,Y,Z\in\mathfrak{X}(M)$, we define
$\mathrm{R}(X,Y)Z=\nabla_{[X,Y]}Z-\nabla_X \nabla_Y Z + \nabla_Y
\nabla_X Z$. The Ricci tensor and the scalar curvature
will be denoted by $\mathrm{Ric}$ and $S$.
The mean curvature $H$ of the hypersurface $M$ is defined by
\[
H=\frac{\alpha}{m}\mathrm{tr}(A)=\frac{1}{m}\sum_{i=1}^{m}\II(E_{i},E_{i})\varepsilon_i\,.
\]
The $(M,g)$-sectional curvature of the plane spanned by two vectors $v_p$ and $w_p$
in $\mathrm{T}_p M$, will be denoted by $K(v_p,w_p)$. The symbols $K^{\II}(v_p,w_p)$
and $\overline{K}(v_p,w_p)$ will be used in accordance with the remark of \S\,\ref{remarkbarIInotation}.
Similarly, the scalar curvature of the second fundamental form will be denoted by $S_{\II}$.

\subsection{The difference tensor $L$.}
The difference tensor $L$ between the two Levi-Civita connections
$\nabla^{\II}$ and $\nabla$ is defined by
\[
L(X,Y) = \nabla^{\II}_{X}Y - \nabla_{X}Y\,,
\]
where $X,Y\in\mathfrak{X}(M)$.
The trace of $L$ with respect to
$\II$ is defined to be the vector field
\[
\mathrm{tr}_{\II}L = \sum_{i=1}^{m}L(V_{i},V_{i})\kappa_i\,,
\]
where $V_i$ and $\kappa_i$ have been defined in \S\,\ref{sec:frame_fields}.

\begin{figure}
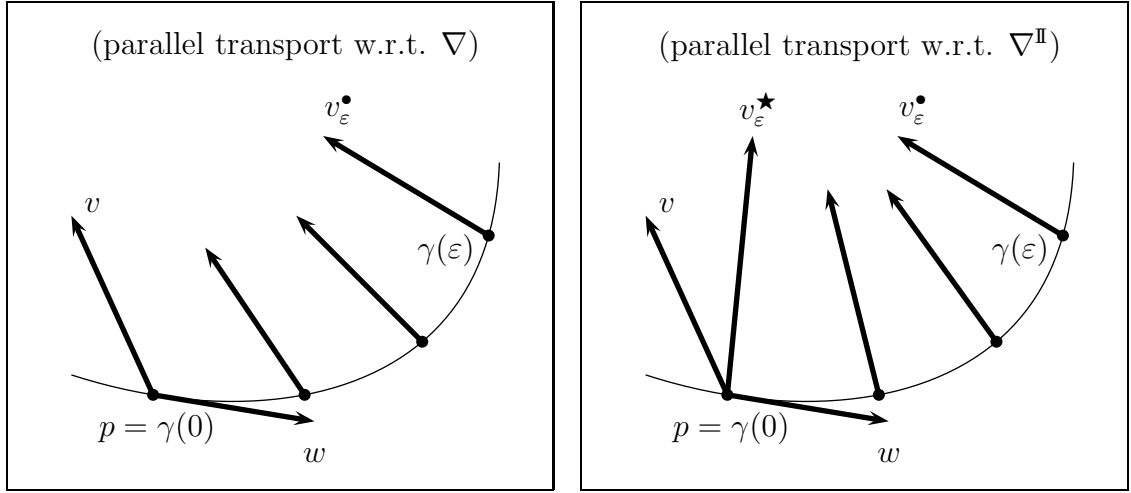

\begin{center}
\begin{tabular}{cc}
%%%begin left figure
\makebox[0.45\textwidth]{\framebox[0.45\textwidth]{
    \rule{0pt}{0.39\textwidth}  %%to ensure the frame will be high enough
    \psset{unit=0.044\textwidth}

%%%curve gamma
\parametricplot[linecolor=black,linewidth=0.5pt]{0}{8}
   {
   t 11 mul sin 8 mul -4 add %%x-coord
   t t mul 0.125 mul -0.5 t mul add 2 add %%y-coord
   }

%%%point p
   \psdot[dotsize=4.5pt](-2.4735,1.625)
    \rput(-2.4,1){$p=\gamma(0)$}
    \psline[linewidth=2pt]{->}(-2.4735,1.625)(-4,5)
    \rput(-3.6,5.2){$v$}
    \psline[linewidth=2pt]{->}(-2.4735,1.625)(0.5417,1.125)
    \rput(0.55,0.5){$w$}

%%first intermediate point
   \psdot[dotsize=4.5pt](.35711, 1.625)
   \psline[linewidth=2pt]{->}(.35711, 1.625)(-1.5,4.4)

%%second intermediate point
   \psdot[dotsize=4.5pt](2.5532, 2.625 )
   \psline[linewidth=2pt]{->}(2.5532, 2.625 )(0.2,5)

%%point gamma(epsilon)
   \psdot[dotsize=4.5pt](3.7949, 4.625 )
   \rput(3,4.35){$\gamma(\varepsilon)$}
   \psline[linewidth=2pt]{->}(3.7949, 4.625 )(0.7,6.5)
   \rput(1,7){$v^{\bullet}_{\varepsilon}$}

\rput(0,8.2){(parallel transport w.r.t. $\nabla$)}

%%to ensure the frame will be high enough
    \rule{0pt}{0.39\textwidth}

}}
%%%end left figure
&
%%%begiun right figure
\makebox[0.45\textwidth]{\framebox[0.45\textwidth]{
    \rule{0pt}{0.39\textwidth}  %%to ensure the frame will be high enough
    \psset{unit=0.044\textwidth}

%%%curve gamma
\parametricplot[linecolor=black,linewidth=0.5pt]{0}{8}
   {
   t 11 mul sin 8 mul -4 add %%x-coord
   t t mul 0.125 mul -0.5 t mul add 2 add %%y-coord
   }

%%%point p
   \psdot[dotsize=4.5pt](-2.4735,1.625)
    \rput(-2.4,1){$p=\gamma(0)$}
    \psline[linewidth=2pt]{->}(-2.4735,1.625)(-4,5)
    \rput(-3.6,5.2){$v$}
    \psline[linewidth=2pt]{->}(-2.4735,1.625)(0.5417,1.125)
    \rput(0.55,0.5){$w$}
    \psline[linewidth=2pt]{->}(-2.4735,1.625)(-2,6.5)
    \rput(-1.9,7){$v^{\bigstar}_{\varepsilon}$}

%%first intermediate point
   \psdot[dotsize=4.5pt](.35711, 1.625)
   \psline[linewidth=2pt]{->}(.35711, 1.625)(-0.6,5.5)

%%second intermediate point
   \psdot[dotsize=4.5pt](2.5532, 2.625 )
   \psline[linewidth=2pt]{->}(2.5532, 2.625 )(0.5,5.5)

%%point gamma(epsilon)
   \psdot[dotsize=4.5pt](3.7949, 4.625 )
   \rput(3,4.35){$\gamma(\varepsilon)$}
   \psline[linewidth=2pt]{->}(3.7949, 4.625 )(0.7,6.5)
   \rput(1,7){$v^{\bullet}_{\varepsilon}$}

\rput(0,8.2){(parallel transport w.r.t. $\nabla^{\II}$)}

%%to ensure the frame will be high enough
    \rule{0pt}{0.39\textwidth}

}}
%%%end right figure
\\
\end{tabular}
\caption{Interpretation of the difference tensor in terms of parallel transport.}
\label{fig:difftensor}
\end{center}
\end{figure}

\begin{remark}
The difference tensor $L$ can easily be interpreted in terms of
parallel transport. Assume $p\in M$ and $v,w\in\mathrm{T}_p M$ are
given. Choose a curve $\gamma:\mathbb{R}\rightarrow M$ such that
$\gamma(0)=p$ and $\gamma'(0)=w$. By $v^{\bullet}_\varepsilon$ we
will denote the vector of $\mathrm{T}_{\gamma(\varepsilon)}M$
obtained by parallel translation of $v$ along $\gamma$ with
respect to $\nabla$. By $v^{\bigstar}_\varepsilon$ we will denote
the vector of $\mathrm{T}_p M$ which is obtained by parallel
transport of the vector $v^{\bullet}_{\varepsilon}$ back to $p$
along $\gamma$ with respect to $\nabla^{\II}$ (see
Figure~\ref{fig:difftensor}). It is not hard to show that
\[
L(v,w)=\lim_{\varepsilon\rightarrow
0}\frac{v^{\bigstar}_{\varepsilon}-v}{\varepsilon}\,.
\]
\end{remark}

%%%begin figure : difference tensor
%%%end figure : difference tensor

\subsection{The equations of Gauss and Codazzi.}
The Riemann-Christoffel curvature tensor
$\mathrm{R}$ of the hypersurface $M$ is related to the second fundamental form
by means of the Gauss equation
\[
g(\mathrm{R}(X,Y)Z,W) = \bg(\overline{\mathrm{R}}(X,Y)Z,W) +
\alpha \left\lgroup \II(X,Z)\,\II(Y,W) - \II(X,W)\,\II(Y,Z) \right\rgroup\,,
\]
which is valid for all tangent vector fields $X,Y,Z,W\in\mathfrak{X}(M)$.
As a consequence, we have
\begin{equation}
\mathrm{Ric}(X,Y) =
\overline{\mathrm{Ric}}(X,Y)-\alpha\, \overline{g}(\overline{\mathrm{R}}{(X,U)}Y,U)
+\alpha\, m \,H \,\II(X,Y)
-\alpha\, \III(X,Y)\,.
\label{eq:gausscontracted}
\end{equation}
The Codazzi equation of the hypersurface is
\[
(\nabla_{X}A)Y - (\nabla_{Y}A)X = \overline{\mathrm{R}}(X,Y)U\ ,
\]
for all $X,Y\in\mathfrak{X}(M)$.

%%%%%%%%%%%%%%%%%%%%%%%%%%%%%%%%%%%%%%%%%%%%%%%%%%%%%%%%%%%%%%%%%%%%%%%%%%%%%%%%%%%%%%%%%%%%%

\section{The Variation of the Area of the Second Fundamental Form}
\label{sec:IImin}

\subsection{The area functional of the second fundamental form.} Let $\E$ denote the set
of all hypersurfaces in a se\-mi-Rie\-mann\-ian manifold $(\bm,\bg)$ for which the first as well as the second fundamental
form is a se\-mi-Rie\-mann\-ian metrical tensor.
Our first objective is to determine the critical points of the area functional of the second fundamental form
\[
\area_{\II} : \E \rightarrow \mathbb{R} : M\mapsto
\area_{\II}(M)=\int_M\,\dd\Omega_{\II}\,.
\]

\subsection{The mean curvature of the second fundamental form.}

\begin{definition}
Let $M$ be a hypersurface in a se\-mi-Rie\-mann\-ian manifold $(\bm,\bg)$, and suppose that
the first as well as the second fundamental form of $M$ is a se\-mi-Rie\-mann\-ian metrical tensor. Let
\[
\mu : \left]-\varepsilon,\varepsilon\right[ \times M \rightarrow \bm : (s,p) \mapsto \mu_s(p)
\]
be a mapping such that
\[
\left\{
\begin{array}{l}
\textrm{$\mu_s(M)\in\E$ for all $s$;}\\
\textrm{$\mu_s(p)=p$ for all $p$ outside of a compact set of $M$ and all $s$;}\\
\textrm{$\mu_0(p)=p$ for all $p\in M$.}
\end{array}
\right.
\]
Then $\mu$ will be called a \textit{variation of $M$ in $\E$}.
\end{definition}

\begin{definition} Let $M$ be a se\-mi-Rie\-mann\-ian
hypersurface of a semi-Riemannian manifold $(\bm,\bg)$ which belongs to the class $\E$.
The vector field $\mathcal{Z}$ in $\mathfrak{X}(M)$ is defined by
\[
\mathcal{Z}= \sum_{i=1}^m \kappa_i \,A^{\leftarrow} \left(\left[\,\overline{\mathrm{R}}(V_i,U)V_i \right]^{T}\right)\,.
\]
Here $A^{\leftarrow}$ denotes the inverse of the shape operator $A$, and $V_i$ and $\kappa_i$ were defined in \S\,\ref{sec:frame_fields}.
\end{definition}

It can easily be seen that the vector field $\mathcal{Z}$ vanishes if $(\bm,\bg)$ has constant sectional curvature.
If $\bm$ has dimension three, the vector field $\mathcal{Z}$ is equal to $\frac{A(Z)}{\dA}$, where the vector field $Z$ has been
defined in \cite{aledohaesenromero,haesenverpoortverstraelen} by the condition
\[
\forall\, X \in\mathfrak{X}(M), \qquad \overline{\textrm{Ric}}(U,X)=\II(Z,X)\,.
\]

\begin{theorem}
Let $M$ be a hypersurface in a se\-mi-Rie\-mann\-ian
manifold $(\bm,\bg)$ for which the first as well as the second fundamental
form is a se\-mi-Rie\-mann\-ian metrical tensor.
Let $\mu$ be a variation of $M$ in $\E$, for which the variational
vector field has normal component $f\, U$.
The variation of the area
functional ${\area}_{\II}$ is given by
\begin{eqnarray*}
\deltasnought {\area}_{\II}(\mu_s{M}) & = & -\alpha\int_{M}
f\cdot\frac{1}{2}\bigg\lgroup
m\, H-\sum_{i=1}^{m}\bg(\overline{\mathrm{R}}(V_i,U)V_i,U)\,\kappa_i
\\
&&\qquad\qquad\qquad
+\frac{\alpha}{2}\Delta_{\II}\mathrm{log}\left|\dA\right|
-\alpha\,\mathrm{div}_{\II}\mathcal{Z}\bigg\rgroup\,
\dd\Omega_{\II}\,.
\end{eqnarray*}
\label{maintheorem}
\end{theorem}
This theorem can be proved by similar methods to those used in \cite{haesenverpoortverstraelen} (see also \cite{verpoort2008}). 
The formula for the variation
of the second fundamental form which was given there, can be generalised to hypersurfaces in the following way:
\[
\deltasnought \II(\mu_s)(X,Y)
 =  \alpha \,f\left\lgroup \bg(\overline{\mathrm{R}}(U,X)U,Y) -\III(X,Y) \right\rgroup +\Hess_{f}(X,Y)\,.
\]
The left-hand side of this expression, which is valid if
the variational vector field is equal to $f U$,
is defined as in \cite{haesenverpoortverstraelen}.

\begin{definition}
\label{def:HII}Let $M$ be an $m$-dimensional hypersurface in a se\-mi-Rie\-mann\-ian manifold $(\bm,\bg)$ for which both the first
and the second fundamental forms are se\-mi-Rie\-mann\-ian metrical tensors. The \textit{mean curvature of the second
fundamental form}
$H_{\II}$ is defined by
\begin{equation}
\label{eq:defhii}
H_{\II} =\frac{1}{2}\left\lgroup
m\, H-\sum_{i=1}^{m}\bg(\overline{\mathrm{R}}(V_i,U)V_i,U)\,\kappa_i
%\right.\\
%&&\qquad\qquad\qquad\left.
+\frac{\alpha}{2}\Delta_{\II}\mathrm{log}\left|\dA\right|
-\alpha\,\mathrm{div}_{\II}\mathcal{Z}\right\rgroup\,.
%\nonumber
\end{equation}
If $H_{\II}=0$, the hypersurface will be called \textit{$\II$-minimal}.
\end{definition}
\begin{remark}
This definition extends those of \cite{glassnersimon,glassner}; in \cite{haesenverpoortverstraelen},
the sign of $H_{\II}$ was chosen differently.
\end{remark}
\begin{example}
The standard embedding of $\mathrm{S}^{m}(\frac{1}{\sqrt{2}})$ in
$\mathrm{S}^{m+1}(1)$ is $\II$-minimal. Furthermore, the standard
embedding of $\mathrm{S}^{k}(\frac{1}{\sqrt{2}})\times\mathrm{S}^{m-k}(\frac{1}{\sqrt{2}})$
in $\mathrm{S}^{m+1}(1)$ (see, e.g., \cite{lawson1969}) is a $\II$-minimal hypersurface
($k=1,\ldots,m-1$). These assertions can be proved with ease when one takes into account the fact that
these hypersurfaces are parallel (in the sense that $\nabla\II=0$).
\end{example}
\begin{remark}
\label{remarkHII}
As a consequence of Theorem~\ref{maintheorem} and Definition~\ref{def:HII}, we obtain the following formulae for the variation of the classical area ($\area$) and
of the area of the second fundamental form ($\area_{\II}$):
\[
\left\{
\begin{array}{rcl}
\displaystyle \deltasnought {\area}(\mu_s(M))        
& \displaystyle= & \displaystyle-m\,\alpha\int f H\,\dd\Omega            \,; \\
 & & \\
\displaystyle\deltasnought {\area}_{\II}(\mu_s(M))  
&\displaystyle = &\displaystyle -\alpha\int f H_{\II}\,\dd\Omega_{\II} \,.
\end{array}
\right.
\]
\end{remark}
\begin{remark}
The expression for $H_{\II}$ can be rewritten in an alternative way at a point $p\in M$ where the
frame fields can be chosen such that
\begin{itemize}
\item
the $g$-orthonormal basis $\{E_1(p),\ldots,E_m(p) \}$
of $\mathrm{T}_p M$ is composed of eigenvectors of the shape operator (principal directions) at $p$:
\[
A(E_i(p))=\lambda_i(p)\,\, E_i(p)\,, \qquad (i=1,\ldots m)\,;
\]
\item
the $\II$-orthonormal basis $\{V_1(p),\ldots,V_m(p) \}$
of $\mathrm{T}_p M$ consists of the rescaled principal directions at $p$:
\[
V_i(p) =  \frac{1}{\sqrt{|\lambda_i(p)|}}E_i(p)\,, \qquad (i=1,\ldots m)\,.
\]
\end{itemize}
In this case, the following expression for the mean curvature of the second fundamental form holds at the point $p$:
\begin{equation}
\label{eq:defhiibis}
\left(H_{\II} \right)_{(p)} = \left( \frac{1}{2}\left\lgroup m H -\sum_{i=1}^{m}\frac{1}{\lambda_i}\overline{K}(E_i,U)\right\rgroup
+\frac{\alpha}{4}\Delta_{\II}\mathrm{log}\left|\dA\right|
-\frac{\alpha}{2}\mathrm{div}_{\II}\mathcal{Z}\right)_{(p)}\,.
\end{equation}
\end{remark}
\begin{remark}
By using the contracted Gauss equation~(\ref{eq:gausscontracted}),
yet another expression for the mean curvature of the second fundamental form can be derived:
\begin{equation}
\label{eq:defhiitris}
H_{\II} =-\frac{\alpha}{2}\Big\lgroup\mathrm{tr}_{\II}\overline{\mathrm{Ric}}
- \mathrm{tr}_{\II}\mathrm{Ric}+\alpha(m^2-2m)H\\
-\frac{1}{2}\Delta_{\II}\mathrm{log}\left|\dA\right|
+\mathrm{div}_{\II}\mathcal{Z}\Big\rgroup\,.
\end{equation}
\end{remark}

%%%%%%%%%%%%%%%%%%%%%%%%%%%%%%%%%%%%%%%%%%%%%%%%%%%%%%%%%%%%%%%%%%%%%%%%%%%%%%%%%%%%%%%%%%%%%

\section{A comparison Result for the Connections}
\label{sec:somefirstresults}

In the sequel of this article we will make use of the following Lemma, which slightly extends well-known results (\cite{hicks1965} Thm. 7,
\cite{simon1972}, and \cite{gardner1972}, Cor. 13). First we recall a useful definition.

\begin{definition}
A totally umbilical, compact  hypersurface $M$ of a se\-mi-Rie\-mann\-ian manifold $(\bm,\bg)$ which
satisfies $A=\rho\,\mathrm{id}$ for a constant $\rho\in\mathbb{R}$, is called an \textit{extrinsic hypersphere}.
\end{definition}
\begin{lemma}
\label{lemma:hicksextension} Let $M$ be a compact
hypersurface of a se\-mi-Rie\-mann\-ian manifold $(\bm,\bg)$. Suppose
that both the first and the second fundamental forms are positive
definite and that these metrical tensors induce the same
Levi-Civita connection. Furthermore, assume that $(M,g)$ has either strictly
positive or strictly negative sectional curvature. Then $M$ is an extrinsic
hypersphere.
\end{lemma}
\begin{proof}\textit{(First Version.)\,}
As an immediate consequence of $\nabla=\nabla^{\II}$, we see that
\[
\mathrm{R}(X,Y)Z=\mathrm{R}^{\II}(X,Y)Z
\]
holds for all $X,Y,Z\in\mathfrak{X}(M)$. Let $p\in M$ be an arbitrary point and choose an orthonormal
basis $\{E_1(p),\ldots,E_m(p)\}$ as in Remark~\ref{remarkHII}:
\[
A(E_i(p))=\lambda_i(p)\, E_i(p) \qquad\qquad (i=1,\ldots, m)\,.
\]
These vectors can be extended to a smooth orthonormal
frame field $\{E_1,\ldots,E_m\}$ on a neighbourhood of $p$ in $M$. For any choice of $i\neq j\in\{1,\ldots,m\}$, there holds
\begin{eqnarray*}
K^{\II}(E_i(p),E_j(p))&=&\left(\frac{\II(\mathrm{R}^{\II}(E_i,E_j)E_i,E_j)}{\II(E_i,E_i)\II(E_j,E_j)}\right)_{(p)}\\
&=& \left(\frac{\alpha\lambda_j g(\mathrm{R}(E_i,E_j)E_i,E_j)}{\lambda_i \lambda_j}\right)_{(p)}\\
&=& \frac{\alpha}{\lambda_i(p)} K(E_i(p),E_j(p))\,.
\end{eqnarray*}
Since the above equation remains valid if the r\^{o}le of $i$ and
$j$ is interchanged whereas $K(E_i(p),E_j(p))\neq0$, it follows that $M$
is totally umbilical. This means that $A=\rho\, \mathrm{id}$ for a
function $\rho:M\rightarrow \mathbb{R}$. Furthermore, for all
$X,Y,Z\in\mathfrak{X}(M)$,
\[
0= \left(\nabla^{\II}_X \II\right)(Y,Z)= \left(\nabla_X \II\right)(Y,Z)=\alpha X\left[\rho\right] g(Y,Z)\,.
\]
Consequently, $\rho$ is a constant.
\end{proof}
\begin{proof} \textit{(Second Version.)\,}
Assume that the conditions, as stated in the lemma, are satisfied for a hypersurface which is not an extrinsic hypersphere. 
Thus the open set of all non-umbilical points on $M$ is non-empty and a non-umbilical point $p\in M$ can be chosen, in a neighbourhood of which the principal directions $\{E_1,\ldots,E_m\}$ are smooth. The fact that $p$ is non-umbilical means, after a possible re-numbering of indices, that
\[
\left\{
\begin{array}{rcl}
\lambda_1(p) &=& \lambda_2 (p) \,= \cdots =\, \lambda_k(p)\,; \\
\lambda_1(p) &\neq& \lambda_{k+1}(p)\,; \\
&\vdots& \\
\lambda_1(p) &\neq& \lambda_{m}(p)\,,
\end{array}
\right.
\]
for some $k\in \left\{1,\ldots,m-1\right\}$\,.

Since the first and the second fundamental form induce the same Levi-Civita connection, there holds $\nabla A=0$. Thus for $X\in\mathfrak{X}(M)$ and $i,j\in\{1,\ldots,m\}$, we have
\[
0 
= 
g \left( (\nabla_X A)E_i,E_j \right)
=
(\lambda_i-\lambda_j)g(\nabla_X E_i, E_j) - X[\lambda_i] \delta_{i\,j}\,,
\]
and in particular, 
\begin{equation}
\label{eq:ch4_nablaXi}
0= g(\nabla_X E_i,E_j)
\qquad
\textrm{(for $1\leqslant i \leqslant k < k+1 \leqslant j \leqslant m$).}
\end{equation}

Now choose a curve $\gamma$ on $M$ satisfying $\gamma(0)=p=\gamma(1)$, and choose a vector $v_p\in\textrm{span}\left\{E_{k+1}(p),\ldots,E_{m}(p)\right\}$. The vector field which is obtained by parallel translation of this vector $v_p$ along the curve $\gamma$ will be denoted by $V=\xi^1 E_1+\cdots+\xi^m E_m$.
Taking the fact (\ref{eq:ch4_nablaXi}) into account, we deduce that
\[
\hspace{-1mm}
\left\{
\begin{array}{rcccccccccl}
0 &\!\!\!=&\!\! g(\nabla_{\gamma'} V, e_1 ) 
&\!\!=&\!\! \left(\xi^1\right)' &\!\!+&\!\! \xi^2\, g(\nabla_{\gamma'}E_2,E_1) 
&\!\!+&\!\! \cdots &\!\!+&\!\! \xi^k\, g(\nabla_{\gamma'}E_k,E_1) \,;\\
&\!\!\vdots&\!\!&\!\!&\!\!&\!\!&\!\!&\!\!&\!\!&\!\!&\!\!\\
0 &\!\!=&\!\! g(\nabla_{\gamma'} V, E_k ) 
&\!\!=&\!\! \left(\xi\right)' &\!\!+&\!\! \xi^1\, g(\nabla_{\gamma'}E_1,E_k) 
&\!\!+&\!\! \cdots &\!\!+&\!\! \xi^{k-1}\, g(\nabla_{\gamma'}E_{k-1},E_k)\,.\\
\end{array}
\right.
\]
The crucial remark which has to be made is that the functions $\xi^{k+1},\ldots,\xi^{m}$ do not enter in these equations because of (\ref{eq:ch4_nablaXi}), whence a system of first-order differential equations in the 
functions $\xi^1,\ldots,\xi^k$ with initial conditions $\xi^1(0)=0,\ldots,\xi^k(0)=0$ is obtained. The conclusion is that $\xi^1\equiv\cdots\equiv\xi^k\equiv0$.

It has thus been shown that the action of the local holonomy group leaves
the subspace $\mathrm{span}\left\{E_{k+1}(p),\ldots,E_m(p)\right\}$
of $\T_p M$ invariant.

Now the local de Rahm theorem can be applied, which gives us that a neighbourhood of $p\in (M,g)$ is isometric to a Cartesian product of Riemannian manifolds. This is in contradiction with the assumption that $(M,g) $ has non-vanishing sectional curvature.
\end{proof}

%%%%%%%%%%%%%%%%%%%%%%%%%%%%%%%%%%%%%%%%%%%%%%%%%%%%%%%%%%%%%%%%%%%%%%%%%%%%%%%%%%%%%%%%%%%%%

\section{Hypersurfaces in a Space Form}
\label{sec:spaceform}

We shall use $\overline{M}^{m+1}_0(\overline{C})$ to denote the following Riemannian manifolds of dimension $m+1$:
\[
\left\{
\begin{array}{lll}
\textrm{the Euclidean hypersphere} & \mathrm{S}^{m+1}(\frac{1}{\sqrt{\overline{C}}}) & \textrm{(for $\overline{C}>0$)\,;}\\
\textrm{the Euclidean space}       & \mathbb{E}^{m+1} & \textrm{(for $\overline{C}=0$)\,;}\\
\textrm{the hyperbolic space}      & \mathrm{H}^{m+1}(\frac{1}{\sqrt{-\overline{C}}}) & \textrm{(for $\overline{C}<0$)\,.}\\
\end{array}
\right.
\]
We shall use $\overline{M}^{m+1}_1(\overline{C})$ to denote the following Lorentzian manifolds of dimension $m+1$:
\[
\left\{
\begin{array}{lll}
\textrm{the de Sitter space}        & \mathrm{S}^{m+1}_1(\frac{1}{\sqrt{\overline{C}}}) & \textrm{(for $\overline{C}>0$)\,;}\\
\textrm{the Minkowski space}        & \mathbb{E}^{m+1}_1 & \textrm{(for $\overline{C}=0$)\,;}\\
\textrm{the anti-de Sitter space}   & \mathrm{H}^{m+1}_1(\frac{1}{\sqrt{-\overline{C}}}) & \textrm{(for $\overline{C}<0$)\,.}\\
\end{array}
\right.
\]
Each of the above se\-mi-Rie\-mann\-ian manifolds has constant sectional curvature $\overline{C}$.

\begin{lemma}
\label{lemma:spaceform}
Let $M$ be a compact se\-mi-Rie\-mann\-ian hypersurface in a se\-mi-Rie\-mann\-ian manifold $(\bm,\bg)$
of constant sectional curvature $\overline{C}$ and dimension $m+1$ (with $m\geqslant2$).
Assume that the second fundamental form of $M$ is positive definite. The inequality
\begin{equation}
\label{eq:ineq}
S_{\II} \leqslant 2\alpha (m-1)  \left\lgroup H_{\II}+\overline{C}\mathrm{tr}A^{\leftarrow}\right\rgroup
\end{equation}
is satisfied if and only if the Levi-Civita connections of the first and the second fundamental forms coincide.
\end{lemma}
\begin{proof}
The following expressions are valid for the curvatures which are involved in the above inequality:
\[
\left\{
\begin{array}{ccl}
\displaystyle 
H_{\II} 
&
\displaystyle
=
&
\displaystyle
\frac{1}{2}\left\lgroup\alpha\,\mathrm{tr}A-\overline{C}\,\mathrm{tr}A^{\leftarrow} \right\rgroup
+\frac{\alpha}{4}\frac{\Delta_{\II}\dA}{\dA}
-\frac{\alpha}{4} \frac{\II(\nabla^{\II}\dA,\nabla^{\II}\dA)}{(\dA)^{2}}\,;\\
& &\\
\displaystyle
S_{\II} 
&
\displaystyle
=
&
\displaystyle 
\alpha(m-1)\left\lgroup\alpha\,\mathrm{tr}A+\overline{C}\,\mathrm{tr}A^{\leftarrow} \right\rgroup+\II(L,L)
-\frac{1}{4} \frac{\II(\nabla^{\II}\dA,\nabla^{\II}\dA)}{(\dA)^{2}}\,,
\end{array}
\right.
\]
where the quantity $\II(L,L)$ is defined by
\[
\II(L,L)=\sum_{i,\,j,\,k\,=1}^{m} (\II(L(V_i,V_j),V_k))^2 \kappa_i\kappa_j\kappa_k =\sum_{i,\,j,\,k\,=1}^{m} (\II(L(V_i,V_j),V_k))^2 \,.
\]
The first expression is an immediate consequence of Equation~(\ref{eq:defhiibis}). The second expression can be
found in, e.g., \cite{schneider1972} (if $(\bm,\bg)$ is the Euclidean space of dimension $m+1$), \cite{aledoromero2003} (if $(\bm,\bg)$ is
the de Sitter space of dimension $m+1$),
or \cite{aledoaliasromero2005} (if $(\bm,\bg)$ is a Riemannian space form of dimension $m+1$).
The inequality~(\ref{eq:ineq}) is equivalent to
\[
0\leqslant \frac{(m-1)}{2}\frac{\Delta_{\II}\dA}{\dA}-\frac{(2m-3)}{4}\frac{\II(\nabla^{\II}\dA,\nabla^{\II}\dA)}{(\dA)^2}-\II(L,L)\,,
\]
and this implies
\[
\dA = \textrm{constant}\qquad\textrm{and}\qquad \nabla=\nabla^{\II}\,.
\]
Conversely, if $\nabla=\nabla^{\II}$, it follows that $\nabla\II$ vanishes. Consequently, $\dA$
is a constant and the inequality is satisfied.
\end{proof}

\begin{theorem}
\label{thm:spaceform}
Let $M$ be a compact Riemannian hypersurface in the space form $\overline{M}^{m+1}_{e}(\overline{C})$ (for $m\geqslant 2$).
Assume that the second fundamental form of $M$ is positive definite.
The inequality
\begin{equation}
\label{eq:ineqb}
S_{\II} \leqslant 2\alpha (m-1)  \left\lgroup H_{\II}+\overline{C}\,\mathrm{tr}A^{\leftarrow}\right\rgroup
\end{equation}
is satisfied if and only if $M$ is an extrinsic hypersphere.
\end{theorem}
\begin{proof} Three cases will be treated separately.
\begin{itemize}
\item[1.] \textit{$\overline{M}^{m+1}_{e}(\overline{C})$ is a Riemannian space form.} It has already
been shown that inequality~(\ref{eq:ineqb}) implies that $M$ is parallel, in the sense that $\nabla\II$ vanishes.
Such hypersurfaces were classified in Theorem 4 of \cite{lawson1969}. If $\overline{C}\geqslant0$, the only
hypersurfaces with a positive definite second fundamental form 
which appear in this classification are the
extrinsic hyperspheres. If $\overline{C}< 0$, the extrinsic hyperspheres are the only compact hypersurfaces in the classification.
\item[2.] \textit{$\overline{M}^{m+1}_{e}(\overline{C})$ is a Lorentzian space form with $\overline{C}\leqslant 0$.} It follows from the Gauss equation that $(M,g)$ has strictly
negative sectional curvature. The result follows from Lemmas~\ref{lemma:hicksextension} and~\ref{lemma:spaceform}.
\item[3.] \textit{$\overline{M}^{m+1}_{e}(\overline{C})$ is the de Sitter space.}
It follows from~(\ref{eq:ineqb}) that $\nabla A$ vanishes. Consequently, $M$ has constant mean curvature and an application of
Theorem~4 of \cite{montiel1988} concludes the proof.
\end{itemize}
\end{proof}

%%%%%%%%%%%%%%%%%%%%%%%%%%%%%%%%%%%%%%%%%%%%%%%%%%%%%%%%%%%%%%%%%%%%%%%%%%%%%%%%%%%%%%%%%%%%%

\section{Hypersurfaces in an Einstein Space}
\label{sec:einstein}

\begin{theorem}
Let $(\bm,\bg)$ be a Riemannian Einstein manifold of dimension $m+1$ (with $m \geqslant 3$) with strictly positive scalar curvature $\overline{S}$.
Any compact hypersurface $M\subseteq\bm$ with positive definite second fundamental form satisfies
\begin{equation}
\label{eq:vwHii}
H_{\II}+m\, \sqrt{\left(\frac{m-2}{m+1}\right)\overline{S}} \geqslant \frac{1}{2}\mathrm{tr}_{\II}\mathrm{Ric}
\end{equation}
if and only if it is an extrinsic hypersphere with 
$ A=\sqrt{ \frac{\overline{S}}{(m-2)(m+1)}}\,\mathrm{id}$. Moreover, in this case there holds 
$ H_{\II}=\sqrt{ \frac{\overline{S}}{(m-2)(m+1)}}$.
\end{theorem}
\begin{proof}
Since $\overline{\mathrm{Ric}}=\frac{\overline{S}}{m+1}\,\overline{g}$, we deduce
that $\mathrm{tr}_{\II}\overline{\mathrm{Ric}}=\frac{\overline{S}}{m+1}\mathrm{tr}A^{\leftarrow}$\,.
Define $\beta$ and $\rho$ by
\[
\beta =\sqrt{\left(\frac{m-2}{m+1}\right)\overline{S}}
\qquad \qquad
\textrm{and}
\qquad \qquad
\rho  =\sqrt{ \frac{\overline{S}}{(m-2)(m+1)}}\,.
\]
Furthermore, the principal curvatures will be denoted by $\lambda_i$ ($i=1,\ldots,m$).
It follows now from~(\ref{eq:defhiitris}) and the assumption~(\ref{eq:vwHii}) that
\begin{eqnarray*}
\int \mathrm{tr}_{\II}\mathrm{Ric}\,\dd\Omega_{\II}
&=&
\int \left\lgroup 2 \,H_{\II}+\beta\sum_{i=1}^{m}\left(\frac{\rho}{\lambda_i}+\frac{\lambda_i}{\rho}\right) \right\rgroup\,\dd\Omega_{\II}\\
&\geqslant& \int 2 \left\lgroup H_{\II}+m\,\beta \right\rgroup \,\dd\Omega_{\II} \geqslant \int \mathrm{tr}_{\II}\mathrm{Ric}\,\dd\Omega_{\II}\,.\end{eqnarray*}
This is only possible if all principal curvatures are equal to $\rho$.
\end{proof}

%%%%%%%%%%%%%%%%%%%%%%%%%%%%%%%%%%%%%%%%%%%%%%%%%%%%%%%%%%%%%%%%%%%%%%%%%%%%%%%%%%%%%%%%%%%%%

\section{Surfaces in a three-dimensional se\-mi-Rie\-mann\-ian Manifold}
\label{sec:3d}

All previous results agree with \cite{haesenverpoortverstraelen} if the surrounding space is three-dimensional
(except for the sign convention of $H_{\II}$).
Moreover, some results can be sharpened. Assume $M\in\E$ and $m=2$. Let $K_{\II}$ denote the
Gaussian curvature of $(M,\II)$. Consequently, the relation $2 K_{\II}=S_{\II}$ is valid.

\begin{theorem}
\label{theo:KII2} Let $M$ be a compact surface in a
three-dimensional se\-mi-Rie\-mann\-ian manifold $(\bm,\bg)$ and suppose
that the first as well as the second fundamental form of $M$ is
positive definite. Suppose that the Gaussian curvature $K$ of $M$
is strictly positive. Then $M$ is an extrinsic hypersphere if
and only if
\begin{equation}
\label{KIIHIIRIC}
K_{\II} \geqslant  \alpha \,H_{\II}+\frac{1}{2}\mathrm{tr}_{\II}\overline{\mathrm{Ric}}\,.
\end{equation}
\end{theorem}
\begin{proof}
Assume first that (\ref{KIIHIIRIC}) is satisfied. A minor adaptation of
the proof of Proposition~5 of \cite{haesenverpoortverstraelen} shows that
$M$ is totally umbilical, and that equality is attained in (\ref{KIIHIIRIC}).
An application of Theorem~6 of \cite{haesenverpoortverstraelen} shows that we have
\[
K_{\II} = \alpha\, H_{\II}+\frac{1}{2}\mathrm{tr}_{\II}\overline{\mathrm{Ric}}-\frac{1}{4}\Delta_{\II}\mathrm{log}(\dA)\,,
\]
and consequently $\dA$ is a constant. The converse follows
since, if $M$ is an extrinsic hypersphere, Theorem~6 of
\cite{haesenverpoortverstraelen} shows that equality holds in
(\ref{KIIHIIRIC}).
\end{proof}
The following Corollary, which follows immediately from the above Theorem and Theorem~\ref{thm:spaceform}, generalises
a result of \cite{manhart1989a,stamou2003}.
\begin{corollary}
Let $M$ be a compact Riemannian surface in the space form
$\overline{M}^{3}_{0}(\overline{C})$ (with $\overline{C}\in\mathbb{R}$) or the de Sitter space.
Assume that the second fundamental form of $M$ is positive definite and that the Gaussian curvature of $(M,g)$ is
strictly positive.
Then either
\[
H_{\II}-\alpha\, K_{\II}+2\frac{\overline{C}H}{K-\overline{C}}
\]
changes sign or $M$ is an extrinsic sphere.
\end{corollary}

%%%%%%%%%%%%%%%%%%%%%%%%%%%%%%%%%%%%%%%%%%%%%%%%%%%%%%%%%%%%%%%%%%%%%%%%%%%%%%%%%%%%%%%%%%%%%

\section{Curves in a se\-mi-Rie\-mann\-ian Surface}
\label{sec:curves}

Let $\gamma:\left]a,b\right[\rightarrow (\bm,\bg):s\mapsto \gamma(s)$ be an arcwise parametrised
time-like or space-like curve in a se\-mi-Rie\-mann\-ian surface. Let $T$ denote the unit tangent vector $\gamma'$ along
$\gamma$.
It will be supposed that $\bg(\bnabla_T T,\bnabla_T T)$ vanishes nowhere.
By virtue of this property, $\gamma$ is sometimes called a Frenet curve.
On the other hand, this requirement precisely means that $\II$ is a se\-mi-Rie\-mann\-ian metrical tensor on $\gamma$.
Let $\{T,U\}$ be the Frenet frame field along $\gamma$:
\[
T=\gamma',\qquad\qquad U=\frac{1}{\sqrt{\left|\bg(\bnabla_T T,\bnabla_T T)\right|}}\bnabla_T T\,.
\]
Further, we set $\beta=\bg(T,T)=\pm1$ and $\alpha=\bg(U,U)=\pm1$.
The geodesic curvature $\kappa$ of $\gamma$ in $(\bm,\bg)$ is determined by the Frenet-Serret formulae:
\begin{equation}
%\label{eq:frenetserret}
\nonumber
\left\lgroup\begin{array}{c} \bnabla_T T \\ \bnabla_T U\end{array}\right\rgroup
=
\left\lgroup\begin{array}{cc} 0 & \beta\kappa \\ -\alpha\kappa & 0 \end{array}\right\rgroup
\left\lgroup\begin{array}{c} T \\ U\end{array}\right\rgroup\,.
\end{equation}
The geodesic curvature $\kappa$ is equal to the
mean curvature of $\gamma\subseteq (\bm,\bg)$. 
The functional which measures the length of a curve with respect to the second fundamental form, which
will be denoted by $\length_{\II}$
instead of $\area_{\II}$, can be computed as the integral
\[
\length_{\II}(\gamma)=\int_{\gamma}\sqrt{\left|\kappa\right|}\,\mathrm{d}s\,.
\]

Let $\overline{K}$ denote the Gaussian curvature of $(\bm,\bg)$.
A calculation shows
\begin{equation}
H_{\II}=\frac{1}{2}\left\lgroup\frac{-\alpha\overline{K}}{\kappa}+\kappa+\frac{\alpha\beta}{4}\left(2\frac{\kappa''}{\kappa^2}-3\frac{(\kappa')^2}{\kappa^3}\right) \right\rgroup\,.
\label{eq:HIIforcurves}
\end{equation}
\begin{example}
A curve $\gamma$ (with $\kappa>0$) in $\mathbb{E}^2$ is
$\II$-minimal if and only if the curvature $\kappa$, when regarded
as a function of the arc-length, satisfies
\begin{equation}
\label{eq:ode1}
4\,\kappa^4 +2\,\kappa\,\kappa''-3(\kappa')^2=0\,.
\end{equation}
This differential equation can most conveniently be solved by introducing the function $\phi= \frac{1}{\kappa}$. With this notation, the above differential equation (\ref{eq:ode1}) is equivalent to
\begin{equation}
\label{eq:ode2}
4-2\,\phi\,\phi''+(\phi')^2=0\,.
\end{equation}
Taking the derivatives of both sides of the above equation (\ref{eq:ode2}), there results $\phi'''=0$. It is now easily verified that the formula
\[
\kappa(s)=\frac{A}{A^2(s+Q)^2+1}\quad\qquad
\textrm{(where $A\in\left]0,+\infty\right[$ and $Q\in\mathbb{R}$)}
\]
describes the general solution of the initial differential equation (\ref{eq:ode1}). This shows that the catenaries are exactly the $\II$-minimal planar curves. (Compare \cite[\S\,27]{blaschke} for the corresponding variational problem for space curves.)
\end{example}
\begin{remark}
It can be asked as well, whether a curve in $\mathbb{E}^2$ can be found which
minimises $\length_{\II}$ among all curves with $\kappa>0$ joining two given points.
This requirement is stronger than merely $\II$-minimality of $\gamma$,
since non-compactly supported fixed-endpoint variations of the curve should also be taken into account.
A simple argument shows that no such minimum exists: if $\gamma_R$ is an arc of a circle of radius $R$ which joins the
two given points, there holds
\[
\lim_{R\rightarrow\infty} \length_{\II}(\gamma_R)= 0\,.
\]
\end{remark}
\begin{example}
For curves on the unit sphere, the equation $H_{\II}=0$ can be
rewritten as
\[
4\kappa^2-4\kappa^4-2\kappa''\kappa+3(\kappa')^2=0\,.
\]
This is equation (4) of \cite{arroyogaraymencia2003}, if the length functional of the second
fundamental form $\length_{\II}$ is interpreted as the so-called curvature energy functional. As is proved and beautifully illustrated
in \cite{arroyogaraymencia2003}, there exists a discrete family of closed, immersed, $\II$-minimal curves on the unit sphere.
Then $S^{1}(\frac{1}{\sqrt{2}})\subseteq S^{2}(1)$ is an embedded ``$\II$-minimal'' curve which belongs to this family. This curve is,
as is remarked in \cite{arroyogaraymencia2003}, actually a local maximum of $\area_{\II}$.
\end{example}

%%%%%%%%%%%%%%%%%%%%%%%%%%%%%%%%%%%%%%%%%%%%%%%%%%%%%%%%%%%%%%%%%%%%%%%%%%%%%%%%%%%%%%%%%%%%%

\section{Geodesic Hyperspheres in a Riemannian Manifold}
\label{sec:geodhyp}

\subsection{Overview of the Results of this Section.}

As a final example we shall investigate the (sufficiently small) geodesic hyperspheres in a Riemannian
manifold, since these provide us with a naturally defined class of hypersurfaces with
a positive definite second fundamental form. We will use the method of power series expansions, hereby relying on computations
of \cite{chenvanhecke,gray1973,grayvanhecke}.

It will first be shown (Theorem \ref{theoremHIIgeodhyp}) that a Riemannian
space for which the value of $H_{\II}$ agrees for every geodesic hypersphere in any of its points
with the corresponding value for a hypersphere of the same radius in a Euclidean space, has to be locally flat.

Let us next denote the geodesic hypersphere of centre $n$ and radius $r$ by
$\mathscr{G}_n(r)$.
It was asked in \cite{grayvanhecke} whether the Riemannian geometry
of the ambient manifold $(\bm,\bg)$ is fully determined by the area functions
\[
\bm\times \left]0,+\infty\right[\rightarrow \mathbb{R} : (n,r) \mapsto \area(\mathscr{G}_n(r))\qquad\qquad\textrm{($r$ sufficiently small)}
\]
of the geodesic hyperspheres. It appears that a decisive answer has not yet been given. Similarly, it may be asked whether a Riemannian manifold for which every
geodesic hypersphere has the same $\II$-area as a Euclidean hypersphere of the same radius,
is locally flat.
In analogy with \cite{grayvanhecke}, we were only able to give an affirmative answer if additional hypotheses
are made (Theorem \ref{thm:areaIIgeodhyp}). For example, the question is answered in the affirmative if the dimension of the ambient manifold does not exceed five.

Let us now explain the technical calculations which result in Theorems~\ref{theoremHIIgeodhyp} and \ref{thm:areaIIgeodhyp}.

\subsection{The co-ordinate system.}
Let $n$ be a point of a Riemannian manifold $(\bm,\bg)$ of dimension $m+1$, and choose a unit vector $e_0\in\mathrm{T}_n \overline{M}$. Consider the geodesic $\gamma$ satisfying $\gamma(0)=n$ and $\gamma'(0)=e_0$. Our purpose is to determine the first few terms in the power series expansion
(in the variable $r>0$) of the value $\left(H_{\II}\right)_{(\gamma(r))}$ which the mean curvature of the second fundamental form of the geodesic hypersphere $\mathscr{G}_n(r)$ of radius $r$ and centre  $n$ assumes in the point $\gamma(r)$. In extension, the letter $r$ will designate also the distance function with respect to the point $n$. It will be assumed throughout that $r>0$ is sufficiently small, in order that everything below is well-defined.

We choose an orthonormal basis
$\{e_0,\ldots,e_m\}$ of $\mathrm{T}_n\overline{M}$ and consider
the associated normal co-ordinate system $\overline{x}=(x^0,\ldots,x^m)$ of $(\bm,\bg)$ centered at $n$:
\[
\overline{x}\left(\overline{\mathrm{exp}}\left(\sum_{s=0}^{m}t^j e_j\right)\right)=(t^0,\ldots,t^m)\,.
\]
For any fixed $r$, a co-ordinate system of $\mathscr{G}_n(r)$ is given by
$x=(x^1,\ldots,x^m)$
in a $\mathscr{G}_n(r)$-neighbourhood of the point $\gamma(r)=\overline{\mathrm{exp}}(r e_0)$.

\begin{figure}
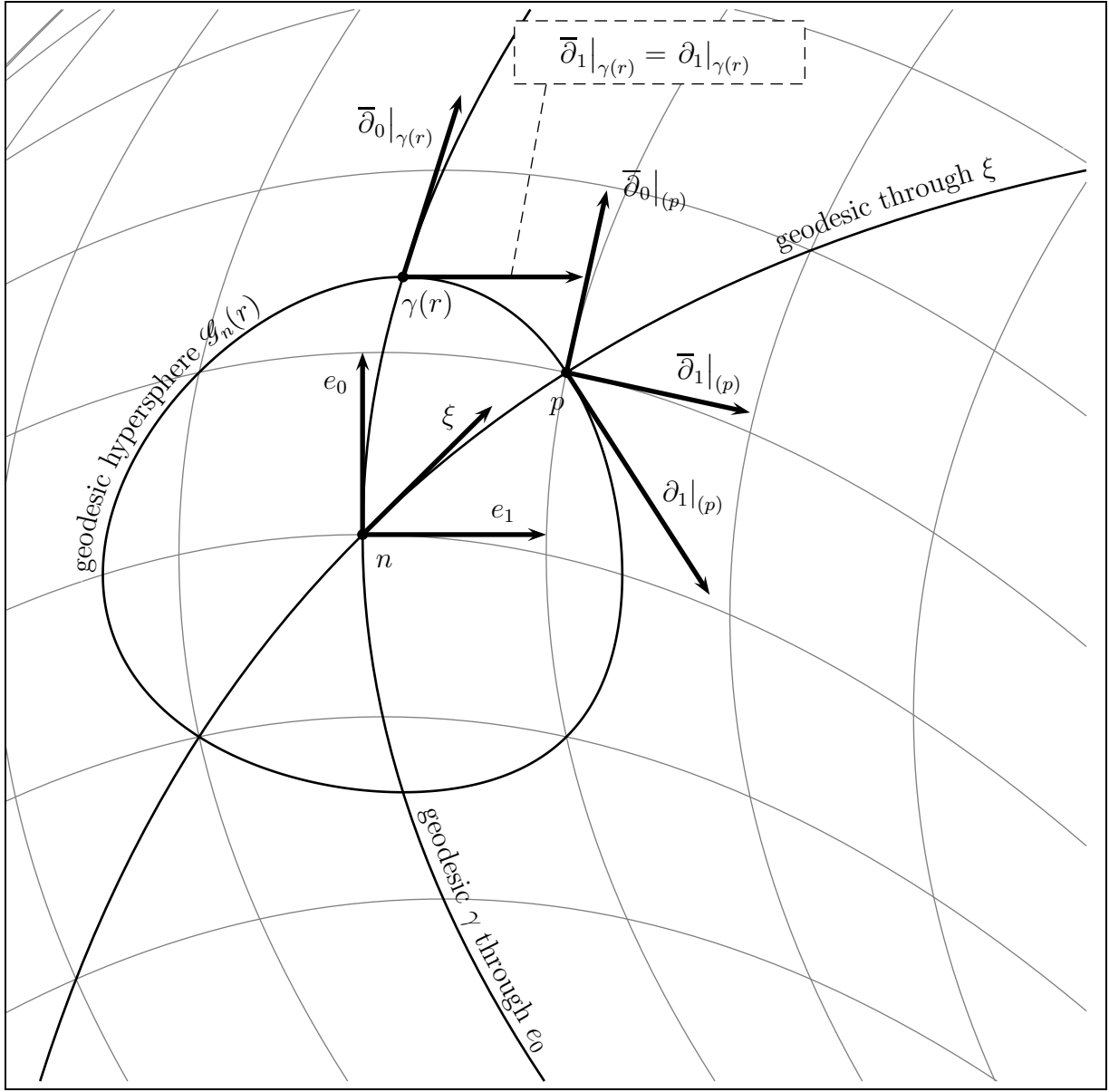

\begin{center}
\makebox[\textwidth]{\framebox[\textwidth]{
    \rule{0pt}{0.98\textwidth}  %%to ensure the frame will be high enough
    \psset{unit=0.098\textwidth}

%%%begin:coordinate lines
%%%coordline x0
\parametricplot[linecolor=gray,linewidth=0.5pt]{-4.6912}{-4.3876}
   {
   2.9750  1.70 t mul add %%x-coord
   13.60 -.1870 t t mul mul add %%y-coord
   }
%%%coordline x1
\parametricplot[linecolor=gray,linewidth=0.5pt]{-3.7012}{-3.1875}
   {
   1.2920 1.70 t mul add %%x-coord
   11.90 -.1870 t t mul mul add %%y-coord
   }
%%%coordline x2 first part
\parametricplot[linecolor=gray,linewidth=0.5pt]{-2.9311}{-1.0342}
   {
   -.0170 1.70 t mul add %%x-coord
   10.20  -.1870 t t mul mul add %%y-coord
   }
%%%coordline x2 second part
\parametricplot[linecolor=gray,linewidth=0.5pt]{1.0342}{2.9512}
   {
   -.0170 1.70 t mul add %%x-coord
   10.20  -.1870 t t mul mul add %%y-coord
   }
%%%coordline x3
\parametricplot[linecolor=gray,linewidth=0.5pt]{-2.3812}{3.5012}
   {
   -.9520 1.70 t mul add %%x-coord
   8.50 -0.1870 t t mul mul add %%y-coord
   }
%%%coordline x4
\parametricplot[linecolor=gray,linewidth=0.5pt]{-2.0512}{3.8312}
   {
   -1.5130 1.70 t mul add %%x-coord
   6.80  -0.1870 t t mul mul add %%y-coord
   }
%%%coordline x5
\parametricplot[linecolor=gray,linewidth=0.5pt]{-1.9412}{3.9412}
   {
   -1.70 1.70 t mul add %%x-coord
   5.10  -0.1870 t t mul mul add %%y-coord
   }
%%%coordline x6
\parametricplot[linecolor=gray,linewidth=0.5pt]{-2.0512}{3.8312}
   {
   -1.5130 1.70 t mul add %%x-coord
   3.40 -.1870 t t mul mul add %%y-coord
   }
%%%coordline x7
\parametricplot[linecolor=gray,linewidth=0.5pt]{-2.3812}{3.0151}
   {
   -.9520 1.70 t mul add %%x-coord
   1.70 -0.1870 t t mul mul add %%y-coord
   }
%%%coordline y1
\parametricplot[linecolor=gray,linewidth=0.5pt]{4.3263}{4.6424}
   {
   -8.50 .1870 t t mul mul add %%x-coord
   2.1080 1.70 t mul add %%y-coord
   }
%%%coordline y2
\parametricplot[linecolor=gray,linewidth=0.5pt]{3.1025}{3.8724}
   {
   -6.80 .1870 t t mul mul add %%x-coord
   3.4170 1.70 t mul add %%y-coord
   }
%%%coordline y3 first part
\parametricplot[linecolor=gray,linewidth=0.5pt]{-2.5600}{-.7313}
   {
   -5.10 .1870 t t mul mul add %%x-coord
   4.3520 1.70 t mul add %%y-coord
   }
%%%coordline y3 second part
\parametricplot[linecolor=gray,linewidth=0.5pt]{.7313}{3.3224}
   {
   -5.10 .1870 t t mul mul add %%x-coord
   4.3520 1.70 t mul add %%y-coord
   }
%%%coordline y4
\parametricplot[linecolor=gray,linewidth=0.5pt]{-2.89}{2.9924}
   {
   -3.40 .1870 t t mul mul add %%x-coord
   4.9130 1.70 t mul add %%y-coord
   }
%%%coordline y5
\parametricplot[linewidth=1pt]{-3}{2.8824}
   {
   -1.70 .1870 t t mul mul add %%x-coord
   5.10 1.70 t mul add %%y-coord
   }
%%%coordline y6
\parametricplot[linecolor=gray,linewidth=0.5pt]{-2.89}{2.9924}
   {
   .1870 t t mul mul %%x-coord
   4.9130 1.70 t mul add %%y-coord
   }
%%%coordline y7
\parametricplot[linecolor=gray,linewidth=0.5pt]{-2.56}{3.3224}
   {
   1.70 .1870 t t mul mul add %%x-coord
   4.3520 1.70 t mul add %%y-coord
   }
%%%coordline y8
\parametricplot[linecolor=gray,linewidth=0.5pt]{-2.01}{2.9251}
   {
   3.40 .1870 t t mul mul add %%x-coord
   3.4170 1.70 t mul add %%y-coord
   }
%%%end:coordinate lines

%%%begin:geodesic circle
\parametricplot[plotpoints=500,linecolor=black,linewidth=1pt]{0}{6.283185}
   {
   -1.70 2.40416 57.29578 t mul sin mul add 57.29578 t mul cos dup mul  0.3740 mul add %%x-coord
    5.10 2.40416 57.29578 t mul cos mul add 57.29578 t mul sin dup mul -0.3740 mul add %%y-coord
   }
%%%end:geodesic circle

%%%begin:geodesic through xi
\parametricplot[linewidth=1pt]{-2.3780}{2.9705}
   {
   -1.70 1.70 t mul add  .1870 t t mul mul add %%x-coord
   5.10  1.70 t mul add -.1870 t t mul mul add %%y-coord
   }
%%%end:geodesic through xi

%%%begin:the point n
    \psdot[dotsize=4.5pt](-1.70,5.10)
    \rput(-1.5,4.85){$n$}

    \psline[linewidth=2pt]{->}(-1.7,5.10)(0,5.1)
    \rput(-0.4,5.3){$e_1$}
          %%e_1
    \psline[linewidth=2pt]{->}(-1.7,5.10)(-1.7,6.8)
    \rput(-1.95,6.5){$e_0$}
          %%e_0
    \psline[linewidth=2pt]{->}(-1.7,5.10)(-0.49792,6.30208)
    \rput(-0.9,6.2){$\xi$}
          %%xi
%%%end:the point n

%%%begin:the point gamma(r)
    \psdot[dotsize=4.5pt](-1.3260,7.5042)
    \rput(-1.1,7.25){$\gamma(r)$}
    \psline[linewidth=2pt]{->}(-1.3260,7.5042)(.345,7.5042)
    \psline[linewidth=2pt]{->}(-1.3260,7.5042)(-.78971,9.2042)
    \rput(-1.4,8.9){$\left.\overline{\partial}_0\right\vert_{\gamma(r)}$}
        %%partial_1
    \psline[linewidth=0.5pt,linestyle=dashed](0,9.3)(-0.3260,7.5042)
    \psframe[fillstyle=solid,fillcolor=white,linecolor=black,linewidth=0.5pt,linestyle=dashed](-0.3,9.3)(2.4,9.9)
    \rput(1,9.55){$\left.\overline{\partial}_1\right\vert_{\gamma(r)}=\left.\partial_1\right\vert_{\gamma(r)}$}
        %%partial_0
%%%end:the point gamma(r)

%%%begin:the point p
    \psdot[dotsize=4.5pt](.1870,6.6130)
    \rput(0.1,6.3){$p$}
    \psline[linewidth=2pt]{->}(.1870,6.6130)(.561,8.313)
    \rput(1,8.313){$\left.\overline{\partial}_0\right\vert_{(p)}$}
        %%overlinepartial_0
    \psline[linewidth=2pt]{->}(.1870,6.6130)(1.8870,6.239)
    \rput(1.5,6.6){$\left.\overline{\partial}_1\right\vert_{(p)}$}
        %%overlinepartial_1
    \psline[linewidth=2pt]{->}(.1870,6.6130)(1.5129,4.5386)
    \rput(1.35,5.45){$\left.\partial_1\right\vert_{(p)}$}
        %%partial_1
%%%end:the point p

%%%begin:text:geodesic hypersphere $G_n(r)$
\pstextpath
   {
   \parametricplot[linecolor=white,linewidth=0pt]{-1.6}{-0.4}
      {
      -1.85  2.54150 57.29578 t mul sin mul add 57.29578 t mul cos dup mul  0.4179 mul add %%x-coord
      5.15   2.54150 57.29578 t mul cos mul add 57.29578 t mul sin dup mul -0.4179 mul add %%y-coord
      }
   }
   {
   geodesic hypersphere $\mathscr{G}_n(r)$
   }
%%%end:text:geodesic hypersphere $G_n(r)$

%%%begin:text:geodesic through xi
\pstextpath
   {
   \parametricplot[linecolor=white,linewidth=0pt]{1.8}{2.9}
      {
      -1.868130000 1.737400000 t mul add  .1870000000 t t mul mul add %%x-coord
      5.268130000  1.737400000 t mul add -0.1870000000 t t mul mul add %%y-coord
      }
   }
   {
   geodesic through $\xi$
   }
%%%end:text:geodesic through xi

%%%begin:text:geodesic through e_0
\pstextpath
   {
   \parametricplot[linecolor=white,linewidth=0pt]{1.45}{3}
      {
      -1.76000000 .1870000000 t t mul mul add %%x-coord
      5.099532500 -1.700000000 t mul add %%y-coord
      }
   }
   {
   geodesic $\gamma$ through $e_0$
   }
%%%end:text:geodesic through e_0

%%to ensure the frame will be high enough
    \rule{0pt}{0.98\textwidth}

}} \caption{A simplified drawing for the co-ordinate systems $x$
and $\overline{x}$. The co-ordinate grid on $(\bm,\bg)$ of $\overline{x}$ is
displayed in gray.} \label{fig:coord}
\end{center}
\end{figure}

It should be noticed that the co-ordinate vector fields $\overline{\partial}_j$ of $\bm$ and $\partial_j$ of $\mathscr{G}_n(r)$
are related by ($j=1,\ldots,m$)
\[
\partial_{j} = \overline{\partial}_{j}-\frac{x^j}{x^0}\overline{\partial}_0\,,
\]
and in particular there holds $\partial_j=\overline{\partial}_j$ along $\gamma$. (See also Figure~\ref{fig:coord}.)
Overlined tensor indices will refer to the co-ordinate system $\overline{x}$, whereas
ordinary tensor indices refer to the co-ordinate system $x$ of the geodesic hyperspheres with
centre $n$.
The coefficients of the
Riemannian curvature tensor of $(\bm,\bg)$ are determined by ($\imath,u,v,e = 0,\ldots,m$)
\[
\R_{\ini\,\inu\,\inv\,\ine}
=
\bg(\R(\overline{\partial}_{\imath},\overline{\partial}_u)\overline{\partial}_{v},\overline{\partial}_{e} )\,.
\]

%%%begin figure : co-ordinate systems
%%%end figure : co-ordinate systems

In \cite{chenvanhecke}, the expansion of the mean curvature $H$
of the geodesic hyperspheres was given at the point $\gamma(r)$:
\begin{eqnarray*}
\left.H\right._{(\gamma(r))} &=&
\frac{1}{r}
-\frac{r}{3m}\left(\Ric_{\innul\,\innul}\right)_{(n)}
-\frac{r^2}{4m}\left(\bnabla_{\innul}\Ric_{\innul\,\innul}\right)_{(n)}\\
&&+\frac{r^3}{m}\left(-\frac{1}{10}\bnabla_{\innul\,\innul}^2\Ric_{\innul\,\innul}-\frac{1}{45}\sum_{a\,e=0}^{m}\left(\R_{\innul\,\ina\,\innul\,\ine}\right)^2\right)_{(n)}
+\Order(r^4)\,.
\end{eqnarray*}
It is follows from this expression that
\textit{the locally flat spaces are the only Riemannian manifolds for which all geodesic hyperspheres have constant mean curvature
which is equal to the inverse of their radius}.

\subsection{The first fundamental form.}
The following expansion for the first fundamental form is given in \cite[Cor. 2.9]{gray1973}:
\begin{eqnarray}
\nonumber
&&\bg_{\ini\,\inj}
=
\delta_{\ini\,\inj}
-\frac{1}{3}\sum_{a\,c=0}^{m} \left(\R_{\ina\,\ini\,\inc\,\inj}\right)_{(n)}  x^a x^c
-\frac{1}{6}\sum_{a\,c\,e=0}^{m}\left(\bnabla_{\ina}\R_{\inc\,\ini\,\ine\,\inj}\right)_{(n)} x^a x^c x^e \\
&&
\label{eq:gij_version1}\quad
+ \frac{1}{120} \sum_{a\,c\,e\,u=0}^{m} \Big(
-6\bnabla_{\ina\,\inc}^{2} \R_{\ine\,\ini\,\inu\,\inj}
+\frac{16}{3}\sum_{s=0}^{m}\R_{\ina\,\ini\,\inc\,\ins}\R_{\ine\,\inj\,\inu\,\ins}
\Big)_{(n)} x^a x^c x^e x^u
+ \Order(r^5)\,.
\end{eqnarray}
This formula is valid for $\imath,\jmath=0,\ldots,m$ and holds on the normal neighbourhood of $n$.
The formula implies
\begin{eqnarray}
\nonumber
\left(\bg_{\ini\,\inj}\right)_{(\gamma(r))}
&=&
\delta_{\ini\,\inj}-\frac{r^2}{3}\left(\R_{\innul\,\ini\,\innul\,\inj}\right)_{(n)}
-\frac{r^3}{6}\left( \bnabla_{\innul}\R_{\innul\,\ini \,\innul \,\inj}\right)_{(n)} \\
&&
\label{eq:gij_version2}
+ \frac{r^4}{120} \left(
-6 \bnabla_{\innul\,\innul}^2 \R_{\innul\, \ini\,\innul\, \inj}
+\frac{16}{3} \sum_{s=0}^{m}\R_{\innul \,\ini\,\innul \,\ins}\R_{\innul \,\inj \,\innul \,\ins}\right)_{(n)}
+ \Order(r^5)\,.
\end{eqnarray}

\subsection{The shape operator of the geodesic hyperspheres.}
It should be noticed that formula~(3.5) of \cite{chenvanhecke} gives the components of the shape operator \textit{with respect to an orthonormal frame field}. As a consequence of this formula~(3.5),
the following expression holds:
\begin{eqnarray}
\label{eq:detA}
&&\left(\log\dA\right)_{(\gamma(r))}+m\log(r)
\,=\,
r^2\left(\frac{-1}{3}\overline{\mathrm{Ric}}_{\innul \, \innul }\right)_{(n)}
+r^3\left(\frac{-1}{4}\bnabla_{\innul }\overline{\mathrm{Ric}}_{\innul \, \innul}\right)_{(n)} \\
\nonumber
&&
\quad+r^4\left(
\frac{-7}{90}\sum_{a\,c=0}^{n}\left(\R_{\innul \, \ina \, \innul \, \inc}\right)^2
-\frac{1}{10}\bnabla_{\innul \,\innul }^2\Ric_{\innul \,\innul }
\right)_{(n)} + \Order(r^5)\,.
\end{eqnarray}
It follows from this equation that \textit{the locally flat spaces are the only Riemannian manifolds for which all geodesic hyperspheres have constant Gauss-Kronecker curvature
which is equal to the inverse of the $m$-th power of their radius.}

In order to find an expression for the \textit{co-ordinate coefficients} of the
shape operator of $\mathscr{G}_n(r)$, we will compute the Christoffel
symbols of $(\bm,\bg)$. Partial derivatives will be denoted with a vertical bar $\vert$ in tensor components.
From (\ref{eq:gij_version1})  we deduce the following expression (for $e,\imath,\jmath=0,\ldots,m$):
\begin{eqnarray}
\nonumber
\left(\bg_{\ini\,\inj\vert\ine }\right)_{(\gamma(r))}
&=&
\frac{-r}{3}\left(\R_{\ine\,\ini\,\innul\,\inj}+
\R_{\innul\,\ini\,\ine\,\inj}\right)_{(n)}\\
&&
\nonumber
-\frac{r^2}{6}\left(\bnabla_{\ine}\R_{\innul\,\ini\,\innul\,\inj}
+ \bnabla_{\innul}\R_{\ine\,\ini\,\innul\,\inj}
+ \bnabla_{\innul}\R_{\innul\,\ini\,\ine\,\inj}
\right)_{(n)}\\
&&
\nonumber
+\frac{r^3}{120}\Bigg(
-6\bnabla_{\ine\,\innul}^2\R_{\innul\,\ini\,\innul\,\inj}
-6\bnabla_{\innul\,\ine}^2\R_{\innul\,\ini\,\innul\,\inj}
-6\bnabla_{\innul\,\innul}^2\R_{\ine\,\ini\,\innul\,\inj}
-6\bnabla_{\innul\,\innul}^2\R_{\innul\,\ini\,\ine\,\inj} \\
&&
\label{eqnarray:partialgij}
\qquad\quad
+\frac{16}{3}\sum_{s=0}^{m}\R_{\ine\,\ini\,\innul\,\ins} \R_{\innul\,\inj\,\innul\,\ins}
+\frac{16}{3}\sum_{s=0}^{m}\R_{\innul\,\ini\,\ine\,\ins} \R_{\innul\,\inj\,\innul\,\ins} \\
&&
\nonumber
\qquad\quad
+\frac{16}{3}\sum_{s=0}^{m}\R_{\innul\,\ini\,\innul\,\ins} \R_{\ine\,\inj\,\innul\,\ins}
+\frac{16}{3}\sum_{s=0}^{m}\R_{\innul\,\ini\,\innul\,\ins} \R_{\innul\,\inj\,\ine\,\ins}
\Bigg)_{(n)}
+\Order(r^4)\,.
\end{eqnarray}
The inverse components of the metric are given by: ($\imath,\jmath=0,\ldots,m$)
\begin{eqnarray}
\label{eqnarray:gij_invers}
\left(\bg^{\ini\,\inj}\right)_{(\gamma(r))}
&=&
\delta_{\ini\,\inj}+r^{2}\left(\frac{1}{3}\R_{\innul\,\ini\,\innul\,\inj}\right)_{(n)}
+r^3\left(\frac{1}{6}\bnabla_{\innul}\R_{\innul\,\ini\,\innul\,\inj}\right)_{(n)}+\Order(r^4)\,.
\end{eqnarray}
\begin{remark}
According to the Gauss lemma, the matrix $(\bg_{\ini\,\inj})$ has the following structure at the point $\gamma(r)$:
\[
(\bg_{\ini\,\inj})_{(\gamma(r))} =
\left\lgroup
\begin{array}{cccc}
1 & 0 & \cdots & 0\\
0 & \bg_{\ineen\,\ineen} & \cdots & \bg_{\ineen\,\inm} \\
0 & \vdots & & \vdots\\
0 & \bg_{\inm\,\ineen} & \cdots & \bg_{\inm\,\inm} \\
\end{array}
\right\rgroup_{(\gamma(r))}
=
\left\lgroup
\begin{array}{cccc}
1 & 0 & \cdots & 0\\
0 & g_{1\,1} & \cdots & g_{1\,m} \\
0 & \vdots & & \vdots\\
0 & g_{m\,1} & \cdots & g_{m\,m} \\
\end{array}
\right\rgroup_{(\gamma(r))}
\,.
\]
Consequently, the same holds for the inverse matrix. This means that (for $\imath,\jmath=1,\ldots,m$) formula~(\ref{eqnarray:gij_invers})
gives also the inverse components
\[
\left(g^{\imath\,\jmath}\right)_{(\gamma(r))}
=
\left(\bg^{\overline{\imath}\,\overline{\jmath}}\right)_{(\gamma(r))}
\]
 of the metrical tensor $g$ of $\mathscr{G}_n(r)$, at a point on the curve $\gamma$.
\end{remark}
The Christoffel symbols $\overline{\Gamma}_{\innul\,\ini}^{\inj}$ of $(\bm,\bg)$ with respect to the
co-ordinate system $\overline{x}$ can be computed
by means of equations~(\ref{eqnarray:partialgij}) and (\ref{eqnarray:gij_invers})
at a point of $\gamma$.

On the other hand, the inward pointing unit normal vector field $U$ of $\mathscr{G}_n(r)$ is given by
\[
U=\frac{-1}{r}\sum_{v=0}^{m} x^v \overline{\partial}_v\,.
\]
Since $\left(r_{\vert\ini}\right)_{(\gamma(r))} = 0$ for $\imath=1\ldots m$, we obtain (for $r > 0$)
\begin{eqnarray*}
A(\left.\partial_{\imath}\right\vert_{(\gamma(r))})
&=&
A(\left.\overline{\partial}_{\imath}\right\vert_{(\gamma(r))})
=
-\left.\bnabla_{\overline{\partial}_{\imath}}(U)\right\vert_{(\gamma(r))}
=
\frac{1}{r} \left.\bnabla_{\overline{\partial}_{\imath}}\left(\sum_{v=0}^{m} x^v \overline{\partial}_v\right)\right\vert_{(\gamma(r))}\\
&=&
\frac{1}{r}\left(\overline{\partial}_{\imath} + \sum_{s\,v=0}^{m}x^v \overline{\Gamma}_{\inv\,\ini}^{\ins}\,\overline{\partial}_s  \right)_{(\gamma(r))}
=
\frac{1}{r}\left(\partial_{\imath} + \sum_{s=0}^{m} r\overline{\Gamma}_{\innul\,\ini}^{\ins}\partial_s\right)_{(\gamma(r))}
\,.
\end{eqnarray*}
Consequently, there holds
$\frac{1}{r}\delta_{i\,s} +\overline{\Gamma}_{\innul\,\ini}^{\ins}= A_{\imath}^{s}$ at the point $\gamma(r)$.
In this way, we obtain the following expression for
the shape operator of $\mathscr{G}_n(r)$ at $\gamma(r)$: ($\imath,\jmath=1,\ldots,m$)
\begin{eqnarray}
\nonumber
\left(A_{\imath}^{s}\right)_{(\gamma(r))}
&=&
\frac{1}{r}\delta_{\imath\,  s}-\frac{r}{3}\left(\R_{\innul\,\ini\,\innul\,\ins} \right)_{(n)}
-\frac{r^2}{4}\left(\bnabla_{\innul}\R_{\innul\,\ini\,\innul\,\ins} \right)_{(n)}\\
&&
\label{eqnarray:Aij}
+r^3 \left(\frac{-1}{10}\bnabla_{\innul\,\innul}^2\R_{\innul\,\ini\,\innul\,\ins}
-\frac{1}{45}\sum_{w=0}^{m}\R_{\innul\,\ini\,\innul\,\inw } \R_{\innul\,\inw\,\innul\,\ins}
\right)_{(n)}
+\Order(r^4)\,.
\end{eqnarray}
Finally, we can compute the components of the second fundamental form in the following way
($\imath,\jmath=1,\ldots,m$):
\begin{eqnarray}
\nonumber
\left(\II_{\imath\,\jmath}\right)_{(\gamma(r))}
&=&
\frac{1}{r}\left(\overline{g}_{\ini\,\inj}\right)_{(n)}- \frac{2r}{3}\left(\R_{\innul\,\ini\,\innul\,\inj}\right)_{(n)}
-\frac{5 r^2}{12}\left( \bnabla_{\innul}\R_{\innul\,\ini\,\innul\,\inj}\right)_{(n)} \\
&&
\label{eqnarray:II_ij}
+r^3\left(\frac{-3}{20}\bnabla_{\innul\,\innul}^2\R_{\innul\,\ini\,\innul\,\inj}
+ \frac{2}{15} \sum_{s=0}^{m}\R_{\innul\,\ini\,\innul\,\ins}\R_{\innul\,\ins\,\innul\,\inj} \right)_{(n)}
+\Order(r^4)\,.
\end{eqnarray}
Since the above equation is only valid at the single point $\gamma(r)=\overline{\mathrm{exp}}(r e_0)$ of $\mathscr{G}_n (r)$,
it needs to be rewritten in order to compute the leading term of $\II_{\imath\,\jmath\vert e}$ at $\gamma(r)$.
A more general expression for $\II_{\imath\,\jmath}$, which is valid at any point $p=\overline{\mathrm{exp}}(r\xi)$
with co-ordinates $(x^0,\ldots,x^m)$ (for a unit vector $\xi\in\mathrm{T}_n \overline{M}$, as in Figure~\ref{fig:coord}), is obtained by
\[
\textrm{substitution of }
\left\{
\begin{array}{l}
\left.\partial_{\imath}\right\vert_{(\gamma(r))}\\
\,\\
\left.\overline{\partial}_{\imath}\right\vert_{(n)}\\
\,\\
e_0
\end{array}
\right.
\textrm{ by }
\left\{
\begin{array}{l}
\left.\partial_{\imath}\right\vert_{(p)}  = \left.\overline{\partial}_{\imath}\right\vert_{(p)}-\frac{x^{\imath}}{x^0}\left.\overline{\partial}_{0}\right\vert_{(p)}\\
\,\\
\left.\overline{\partial}_{\imath}\right\vert_{(n)}-\frac{x^{\imath}}{x^0}\left.\overline{\partial}_{0}\right\vert_{(n)}\\
\,\\
\xi=\frac{1}{r}\sum_{a=0}^{m}x^a e_a\\
\end{array}
\right.
\]
in the previous formula.
The result is
\begin{eqnarray}
\nonumber
\II_{\imath\,\jmath}\!\!\!\!
&=&\!\!\!
\frac{1}{r}\left(
\delta_{\imath\,\jmath}
+\frac{x^{\imath} x^{\jmath} }{(x^0)^2}
-\frac{2}{3}\sum_{a\,c=0}^{m}\left(\R_{\ina\,\ini\,\inc\,\inj} \right)_{(n)} x^a x^c
+\frac{2}{3}\sum_{a\,c=0}^{m}\left(\R_{\ina\,\innul\,\inc\,\inj} \right)_{(n)} 
\frac{x^{\imath}}{x^0}x^a x^c
\right.
\\
\label{eqnarray:II_ij_bis}
&&
\left.\hspace{7mm}
+\,\frac{2}{3}\sum_{a\,c=0}^{m}\left(\R_{\ina\,\ini\,\inc\,\innul} \right)_{(n)} 
\frac{x^{\jmath}}{x^0}x^a x^c
-\frac{2}{3}\sum_{a\,c=0}^{m}\left(\R_{\ina\,\innul\,\inc\,\innul} \right)_{(n)} 
\frac{x^{\imath}\,x^{\jmath}}{(x^0)^2}x^a x^c
\right)
+\Order(r^2)
\,,
\end{eqnarray}
where the function $(x^0)^2$ can be expressed in the co-ordinate system $x$ on $\mathscr{G}_n(r)$ by
\[
(x^0)^2= r^2-(x^1)^2-\ldots-(x^m)^2\,.
\]
Consequently, there holds ($\imath,\jmath,e=1,\ldots,m$):
\begin{equation}
\left(\II_{\imath\,\jmath\vert e}\right)_{(\gamma(r))}
=
\frac{-2}{3}\left(\R_{\ine\,\ini\,\innul\,\inj}+\R_{\innul\,\ini\,\ine\,\inj} \right)_{(n)} + \Order(r)\,.
\end{equation}
In this way, we obtain an expression for the leading term of the Christoffel
symbols of the second fundamental form at $\gamma(r)$ with respect
to the co-ordinate system $x$ of the geodesic hyperspheres ($\imath,\jmath,s=1,\ldots,m$):
\begin{equation}
\label{eq:ChristII}
\left(\left.\Gamma_{\II}\right._{ij}^{s}\right)_{(\gamma(r))}
=
\frac{2 r}{3}\left( \R_{\ins\,\ini\,\innul\,\inj} + \R_{\innul\,\ini\,\ins\,\inj} \right)_{(n)} +\Order(r^2)\,.
\end{equation}
After some work, it can be concluded from equations (\ref{eq:detA}), (\ref{eqnarray:II_ij}) and (\ref{eq:ChristII})
that
\begin{eqnarray*}
\left.\Delta^{\II}\log\dA\right\vert_{(\gamma(r))}
&=&
\frac{-2r}{3}\left(\overline{S}-(m+1)\Ric_{\innul\,\innul}\right)_{(n)}\\
&&+r^2\left(
-\overline{S}_{\vert \innul}
+\frac{3}{4}(m+2)\bnabla_{\innul}\Ric_{\innul\,\innul}
\right)_{(n)}
\\
&&
+r^3\Bigg( \frac{-16}{45}\sum_{v\,w=0}^{m}\R_{\innul\,\inv\,\innul\,\inw}\Ric_{\inv\,\inw}
+\frac{14}{45}(3+m)
\sum_{v\,w=0}^{m}\left(\R_{\innul\,\inv\,\innul\,\inw}\right)^2\\
&& \qquad\qquad- \frac{7}{15}\sum_{\imath\,v\,w=0}^{m}\left(\R_{\ini\,\inv\,\innul\,\inw}\right)^2 -\frac{3}{5}\overline{\mathrm{Hess}}_{(\overline{S})\innul\,\innul}\\
&& \qquad\qquad+ \frac{(6+2m)}{5}\bnabla_{\innul\,\innul}^2 \Ric_{\innul\,\innul}+\frac{22}{45}\sum_{v=0}^{m}\left(\Ric_{\innul\,\inv}\right)^2\\
&& \qquad\qquad -\frac{4}{9}\left(\Ric_{\innul\,\innul}\right)^2-\frac{1}{5}\overline{\Delta}\,\Ric_{\innul\,\innul}\Bigg)_{(n)} +\Order(r^4)\,.
\end{eqnarray*}

\subsection{Further computations}

We will not give the details of the further calculations which can be performed in a similar way. The
$\II$-divergence of the vector field $\mathcal{Z}$
is given by:
\begin{eqnarray*}
\left.\mathrm{div}_{\II}\mathcal{Z}\right\vert_{(\gamma(r))}
&=&  r\left((m+1)\Ric_{\innul\,\innul}-\overline{S}\right)_{(n)}\\
&&+r^2\left((m+2)\bnabla_{\innul}\Ric_{\innul\,\innul}-\frac{3}{2}\overline{S}_{\vert\innul}\right)_{(n)}\\
&&
+r^3\Bigg(
\frac{-1}{3}\sum_{\imath\,\jmath=0}^{n}\R_{\innul\,\ini\,\innul\,\inj}\Ric_{\ini\,\inj}
+\frac{(m+3)}{2}\bnabla_{\innul\,\innul}^{2}\Ric_{\innul\,\innul} \\
&&\qquad\qquad +\frac{2}{3}\sum_{v=1}^m \left(\Ric_{\innul\,\inv}\right)^2
+\frac{(m+3)}{3}\sum_{i\,j=0}^m \left(\R_{\innul\,\ini\,\innul\,\inj}\right)^2 \\
&&\qquad\qquad -\overline{\mathrm{Hess}}_{(\overline{S})\innul\,\innul}
-\frac{1}{2}\sum_{a\,c\,e=0}^{m}\left(\R_{\ina\,\inc\,\ine\,\innul}\right)^2
\Bigg)_{(n)}+\Order(r^4)\,.
\end{eqnarray*}
The $\II$-trace of the $\overline{\mathrm{Ricci}}$ tensor can be calculated as
\begin{eqnarray*}
\left.\mathrm{tr}_{\II}\overline{\mathrm{Ric}}\right\vert_{(\gamma(r))}
&=& r\left(\overline{S}-\Ric_{\innul\,\innul}\right)_{(n)}
+r^2\left(\overline{S}_{\vert\innul}-\bnabla_{\innul}\Ric_{\innul\,\innul}\right)_{(n)}\\
&&\hspace{-2cm} + r^3\left(
\frac{1}{3}\sum_{\imath\,\jmath=0}^{m}\R_{\innul\,\ini\,\innul\,\inj}\Ric_{\ini \,\inj}
-\frac{1}{2}\bnabla^2_{\innul\,\innul}\Ric_{\innul\,\innul}
+\frac{1}{2}\overline{\mathrm{Hess}}_{(\overline{S})\innul\,\innul}
\right)_{(n)}
+\Order(r^4)\,.
\end{eqnarray*}
The $\II$-trace of the Ricci tensor satisfies
\begin{eqnarray*}
\left.\mathrm{tr}_{\II}{\mathrm{Ric}}\right\vert_{(\gamma(r))}
&=& \frac{m(m-1)}{r}+r \left(\overline{S}-\frac{(m+5)}{3}\Ric_{\innul\,\innul}\right)_{(n)}\\
&&
+r^2\left(\overline{S}_{\vert\innul}-\frac{(m+7)}{4}\bnabla_{\innul}\Ric_{\innul\,\innul} \right)_{(n)}\\
&&
+r^3\Bigg(
\frac{1}{3}\sum_{\imath\,\jmath=0}^{m}\R_{\innul\,\ini\,\innul\,\inj}\Ric_{\ini \,\inj}
-\frac{(m+9)}{10}\bnabla^2_{\innul\,\innul}\Ric_{\innul\,\innul}\\
&&\qquad\qquad -\frac{(m+14)}{45}\sum_{\imath\,\jmath=0}^{m} \left(\R_{\innul\,\ini\,\innul\,\inj}\right)^2
+\frac{1}{2}\overline{\mathrm{Hess}}_{(\overline{S})\innul\,\innul}
\Bigg)_{(n)}+\Order(r^4)\,.
\end{eqnarray*}

\subsection{An expression for $H_{\II}$.}
From the previous computations and formula~(\ref{eq:defhiitris}), we obtain
\begin{eqnarray}
\left.H_{\II}\right._{(\gamma(r))} &=& \frac{m}{2r}
+\frac{r}{3}\left(\overline{S}-(m+3)\Ric_{\innul\,\innul}\right)_{(n)}\nonumber\\
&&+r^{2}\left(\frac{1}{2}\overline{S}_{\vert\innul}-\frac{(20+5m)}{16}\bnabla_{\innul}\Ric_{\innul\,\innul}\right)_{(n)}\nonumber\\
&&
+r^3\Bigg(
\frac{7}{90}\sum_{\imath\,\jmath=0}^{m}\R_{\innul\,\ini\,\innul\,\inj}\Ric_{\ini \,\inj}
-\frac{(15+3m)}{20}\bnabla^2_{\innul\,\innul}\Ric_{\innul\,\innul}\nonumber\\
\label{eq:HIIpowerseries}
&& \qquad \qquad -\frac{19}{90}\sum_{v=1}^m \left(\Ric_{\innul\,\inv}\right)^2
-\frac{(20+4m)}{45}\sum_{\imath\,\jmath=0}^{m} \left(\R_{\innul\,\ini\,\innul\,\inj}\right)^2\\
&& \qquad\qquad +\frac{7}{20}\overline{\mathrm{Hess}}_{(\overline{S})\innul\,\innul}
+\frac{2}{15}\sum_{a\,c\,e=0}^{m}\left(\R_{\ina\,\inc\,\ine\,\innul}\right)^2\nonumber\\
&& \qquad \qquad
+\frac{1}{90}\left(\Ric_{\innul\,\innul}\right)^2
-\frac{1}{20}\overline{\Delta}\,\Ric_{\innul\,\innul}\Bigg)_{(n)}+\Order(r^4)\,. \nonumber
\end{eqnarray}
\begin{theorem}
\label{theoremHIIgeodhyp}
A Riemannian manifold (of dimension $m+1$) is locally flat if and only if
the mean curvature of the second fundamental form
of every geodesic hypersphere is equal to the constant $\frac{m}{2r}$ (where $r$ is the radius of the geodesic hypersphere).
\end{theorem}
\begin{proof}
Suppose that $(\bm,\bg)$ is a Riemannian manifold for which the relation $H_{\II}=\frac{m}{2r}$ holds for every geodesic hypersphere.
Then for any choice of $n\in \bm$ and $e_0 \in\mathrm{T}_n \bm$,
the coefficients of the positive powers of $r$ in formula~(\ref{eq:HIIpowerseries})  vanish.
An analysis of the equation
\[
\forall\,n\in\bm\quad\forall \,e_0\in\mathrm{T}_n\bm\quad\textrm{with}\quad \|e_0\|=1,
\quad\quad (m+3)\overline{\textrm{Ric}}(e_0,e_0)=\overline{S}_{(n)}
\]
gives that $\bm$ is Ricci flat.
The fact that the coefficient of $r^3$ vanishes implies that for each point $n\in\bm$ and
for each unit vector $\xi\in\mathrm{T}_n \bm$, there holds
\[
\frac{(20+4m)}{45}\sum_{\imath\,\jmath=0}^{m} \left(\R_{\xi\,\ini\,\xi\,\inj}\right)^2 =
\frac{2}{15}\sum_{a\,c\,e=0}^{m}\left(\R_{\ina\,\inc\,\ine\,\xi}\right)^2\,.
\]
Both sides of the above equation can be integrated over the unit hypersphere of $\mathrm{T}_n \bm$ with the help of the
results of \cite{chenvanhecke,grayvanhecke}. By means of the resulting equation, it can be concluded that
$\overline{\mathrm{R}}$ vanishes.
\end{proof}

\subsection{The area of geodesic hyperspheres, as measured by means of the second fundamental form}

Let $\alpha_m$ denote the area of a unit hypersphere in $\mathbb{E}^{m+1}$.
A calculation gives
\begin{eqnarray}
\area_{\II}(\mathscr{G}_n(r)) &=& r^{\frac{m}{2}}\alpha_m \Bigg[ 1 - r^2 \left(\frac{\overline{S}}{3(m+1)}\right)_{(n)}\nonumber\\
\label{eq:FIIpowerseries}
&&\qquad\qquad +r^4 \frac{1}{(m+1)(m+3)}\Bigg( \frac{1}{18}(\overline{S})^2+\frac{1}{15}\sum_{\imath\,\jmath=0}^{m}\left(\Ric_{\ini\,\inj}\right)^2\\
\nonumber
&&\qquad\qquad\qquad\qquad -\frac{1}{15}\sum_{a\,c\,e\,s=0}^{m}\left(\R_{\ina\,\inc\,\ine\,\ins}\right)^2-\frac{3}{20}\overline{\Delta}\,\overline{S}\Bigg)_{(n)}
+\Order(r^5)\Bigg]\,.
\end{eqnarray}
Thus, the following adaption of theorem 4.1 of \cite{grayvanhecke} can be obtained.
\begin{theorem}
\label{thm:areaIIgeodhyp}
Let $(\bm,\bg)$ be a Riemannian manifold of dimension $m+1$, and suppose that the area of every geodesic hypersphere
of $\bm$, as seen in the geometry of the second fundamental form, is equal to $r^{\frac{m}{2}}\alpha_m$ (where $r$ is the
radius of the geodesic hypersphere, and $\alpha_m$ is the area of the unit hypersphere of $\mathbb{E}^{m+1}$). Then there holds
\begin{equation}
\label{cond}
\left\{
\begin{array}{rcl}
\overline{S}&=&0\,;\\
\|\R\|^2&=&\|\Ric\|^2\,.
\end{array}
\right.
\end{equation}
Further, $\bm$ is locally flat if any of the following additional hypotheses is made:
\begin{itemize}
\item[(i).] $\dim \bm\leqslant 5$;
\item[(ii).] the Ricci tensor of $\bm$ is positive or negative semi-definite (in particular if $M$ is Einstein);
\item[(iii).] $\bm$ is conformally flat and $\dim\bm\neq6$;
\item[(iv).] $\bm$ is a K\"{a}hler manifold of complex dimension $\leqslant 5$;
\item[(v).] $\bm$ is a Bochner flat K\"{a}hler manifold of complex dimension $\neq 6$;
\item[(vi).] $\bm$ is a product of surfaces (with an arbitrary number of factors).
\end{itemize}
\end{theorem}
\begin{proof}
The first part of the theorem follows immediately from the given power series expansion
(\ref{eq:FIIpowerseries}). Now assume (\ref{cond}) is satisfied.
\begin{itemize}
\item[(i)] Suppose that $\bm$ has dimension $\leqslant5$ (i.e., $m\leqslant4$). The trivial case $m=1$ should be excluded in the reasoning below.
Let $\overline{\mathrm{W}}$ denote the Weyl conformal curvature tensor of $(\bm,\bg)$. There holds
\begin{eqnarray*}
0 &\leqslant& \|\overline{\mathrm{W}}\|^2\\
&=& \|\overline{\mathrm{R}}\|^2-\frac{4}{m-1}\|\overline{\mathrm{Ric}}\|^2+\frac{2}{m(m-1)}\overline{S}^2\\
&=& \frac{m-5}{m-1}\|\overline{\mathrm{R}}\|^2 \leqslant 0  \,,
\end{eqnarray*}
and consequently, $0=\overline{\mathrm{R}}$.
\item[(ii)] If $\epsilon\overline{\mathrm{Ric}}$ is positive semi-definite, for , for $\epsilon=\pm1$, then $0\leqslant\epsilon\mathrm{tr}\overline{\mathrm{Ric}}=\epsilon \overline{S}=0$
and consequently $\overline{\mathrm{Ric}}=0$ and $\overline{\mathrm{R}}=0$.
\item[(iii)] The case where $\dim \bm\leqslant 5$ has already been proved. So assume $\bm$ is a conformally flat Riemannian manifold
which satisfies (\ref{cond}),  $\dim\bm\geqslant 7$ (i.e. $m\geqslant6$) and $0\neq \|\R\|$.
The fact that $0=\|\overline{\mathrm{W}}\|^2$ implies
\[
(m-1) \|\R\|^2 = 4\|\Ric\|^2=4\|\R\|^2 < (m-1) \|\R\|^2\,,
\]
which is clearly a contradiction.
\end{itemize}
(iv) and (v) can be proved similarly to the two previous cases by an analysis of the squared norm of the Bochner curvature tensor. (vi) can
be proved in the same way as in \cite{grayvanhecke}.
\end{proof}

\begin{remark}
For a given $r>0$ and $n\in\bm$, the collection concentric geodesic hyperspheres $\left\{\mathscr{G}_n(r+s)\right\}$
can be seen as a variation of $\mathscr{G}_n(r)$ with variational vector field $-U$. An application of Theorem~\ref{maintheorem}
gives that the relation
\[
\frac{\partial}{\partial r}\area_{\II}(\mathscr{G}_{n}(r))= \int_{\mathscr{G}_{n}(r)} H_{\II}\,\dd\Omega_{\II}
\]
holds. It can indeed be checked that the first terms in the power series expansion of both functions agree.
\end{remark}

%%%%%%%%%%%%%%%%%%%%%%%%%%%%%%%%%%%%%%%%%%%%%%%%%%%%%%%%%%%%%%%%%%%%%%%%%%%%%%%%%%%%%%%%%%%%%

\section*{Acknowledgements}

Some of the above results have been obtained with the help of useful comments made by B.-Y. Chen, F. Dillen, L. Vanhecke and L. Verstraelen. 
Furthermore, J. Bolton was so kind to suggest several linguistic improvements. We wish to express our sincerest gratitude to all of them.

%%%%%%%%%%%%%%%%%%%%%%%%%%%%%%%%%%%%%%%%%%%%%%%%%%%%%%%%%%%%%%%%%%%%%%%%%%%%%%%%%%%%%%%%%%%%%

\end{document}